\ifx\selectfont\undefined
\documentstyle[12pt]{article}
\else
\documentstyle[12pt,oldlfont]{article}
\fi

\newcount\hh
\newcount\mm
\mm=\time
\hh=\time
\divide\hh by 60
\divide\mm by 60
\multiply\mm by 60
\mm=-\mm
\advance\mm by \time
\def\hhmm{\number\hh:\ifnum\mm<10{}0\fi\number\mm}
\makeatletter
\@addtoreset{equation}{section}
\ifx\reset@font\undefined
  \let\reset@font=\relax
\fi
\def\section{\@startsection {section}{1}{\z@}{-3.5ex plus-1ex minus
    -.2ex}{2.3ex plus.2ex}{\reset@font\large\bf}}
\def\subsection{\@startsection{subsection}{2}{\z@}{-3.25ex plus-1ex
    minus-.2ex}{-.1em}{\reset@font\large\bf}}
\def\subsubsection{\@startsection{subsubsection}{3}{\z@}{-3.25ex plus
 -1ex minus-.2ex}{-.1em}{\reset@font\normalsize\bf}}
\makeatother

\title{Asymptotic infinite-dimensional theory of Banach spaces}
\author{Bernard Maurey 
         \and Vitali D.~Milman
        \and  Nicole Tomczak-Jaegermann}
\newcommand\address{\noindent\leavevmode%
Equipe d'Analyse et Math\'ematiques Appliqu\'ees,\\
Universit\'e de Marne la Vall\'ee,\\
93166 Noisy Le Grand CEDEX, France,\\
{\small\tt%
   maurey@logique.jussieu.fr}\\[.5cm]
\noindent 
Department of Mathematics,\\
Raymond and Beverly Sackler Faculty of Exact Sciences\\
Tel Aviv University, Tel Aviv, Israel\\
{\em and} \\
Department of Mathematics, Ohio State University,\\
Columbus, Ohio, 43210, USA,\\
{\small\tt%
   vitali@math.tau.ac.il}\\[.5cm]
\noindent
Department of Mathematics,
University of Alberta,\\
Edmonton, Alberta, Canada T6G 2G1,\\
{\small\tt%
  ntomczak@approx.math.ualberta.ca} }

\date{}

\newtheorem{thm}{Theorem}

\newtheorem{lemma}{Lemma}
\newtheorem{cor}{Corollary}

\newtheorem{dfn}{Definition}

\newbox\nrmbox
\setbox\nrmbox=\hbox{$\Vert$}
\def\nrmrule{\vrule height\ht\nrmbox depth1.2\dp\nrmbox}
\setbox\nrmbox=%
  \hbox{\kern0.15em\nrmrule\kern0.15em\nrmrule\kern0.15em\nrmrule\kern0.15em}
\newcommand{\Snorm}[1]%
  {\copy\nrmbox#1\copy\nrmbox\kern-0.03em\lower.4ex\hbox{}}
\newcommand{\proof}{{\noindent\bf Proof{\ \ }}}
\newcommand{\proofsec}{{\bf Proof{\ \ }}}
\newcommand{\qed}{\bigskip\hfill\(\Box\)}
\newcommand{\Rn}[1]{\mbox{{\it I\kern -0.25emR}$^{\,{#1}}$}}
\newcommand{\NN}{\mbox{{\it I\kern -0.25emN}}}

\newcommand{\al}{\alpha}

\newcommand{\ep}{\varepsilon}

\newcommand{\la}{\lambda}

\newcommand{\di}{{\rm d}}
\newcommand{\spn}{{\rm span\,}}
\newcommand{\pS}{\mbox{\bf S}}
\newcommand{\pV}{\mbox{\bf V}}

\newcommand{\supp}{{\rm supp\,}}

\newcommand{\codim}{{\rm codim\,}}

\newcommand{\cB}{{\cal B}}
\newcommand{\ie}{{\em i.e.,}}
\newcommand{\eg}{{\em e.g.,}}
\newcommand{\cf}{{\em cf.}}

\begin{document}

\maketitle


\def\nek{,\ldots,}

In this paper we study structural properties of infinite dimensional
Banach spaces.  The classical understanding of such properties was
developed in the 50s and 60s; goals of the theory had direct roots in
and were natural expansion of problems from the times  of Banach.  Most
of surveys and books of that period directly or indirectly discussed
such problems as the existence of unconditional basic sequences, the
$c_0$-$\ell_1$-reflexive subspace problem and others.
However, it has been realized recently that such a nice and elegant
structural theory does not exist.  Recent examples (or
counter-examples to classical problems) due to Gowers and Maurey [GM]
and  Gowers [G.2], [G.3] showed much more diversity in the structure 
of  infinite dimensional subspaces of Banach spaces than was expected.

On the other hand, structure of finite dimensional subspaces of
Banach spaces and related local properties have been well understood
in the last two decades.  Many exciting theorems on the behaviour of
high dimensional subspaces, finite rank operators, quotient spaces and
others were discovered.  They have an asymptotic nature: dimension
should increase to infinity to reveal regularities behind an
increasing diversity of discussed objects (\cf~\eg~[MiSch], [P.1],
[P.2], [T]).

In this paper infinite dimensional phenomena are investigated by using
a similar asymptotic approach.  To envisage such phenomena, we discard
all informations of a finite dimensional nature and study
properties of a space ``at infinity''.  This naturally  motivates  a
fundamental concept of {\it asymptotic finite-dimensional spaces\/} of $X$,
which will be explained later in this introduction.
The main idea behind it  is a stabilization at   infinity of
finite dimensional subspaces which appear everywhere far away.
This further leads to an infinite-dimensional construction resulting in a
notion of an {\it asymptotic version} of $X$.

Similar stabilization procedures in the form of the notions of
spectrum and tilda-spectrum, were studied back at the end of the 60s
in [Mi1] and [Mi2].
Originally a complete stabilization procedure was used in [Mi.1] for
stabilizing special geometric moduli, so-called $\beta$- and
$\delta$-moduli (see  also  recent applications of this approach in
[MiT]). This was achieved by considering  functions
$f_\lambda (x, y) = \|x + \lambda y\|$
on the unit sphere $ S(X)$ of $X$; 
in the case when a complete stabilization of these functions
on $S(X)$ was possible,  $X$ was shown to contain 
an infinite-dimensional $l_p$
subspace.  Different, although close, asymptotic view was taken in
[KM] through the notion of stable spaces.  Again, it was proved that
a stable space $X$  contains an infinite-dimensional $l_p$
subspace.  In both approaches strong stabilization conditions allowed
a complete recovery of some infinite-dimensional subspaces through a
construction of what we would call today a ``stabilized asymptotic
version''.

We would like to mention in this context that for spaces $l_p$
for $1 < p < \infty$ there exist  Lipschitz functions
on the sphere (in fact, equivalent norms) which
do not allow a complete stabilization on any infinite-dimensional
subspace of $l_p$. This is a weak form of a
recent distortion result
by Odell and Schlumprecht [OS.2], which which gives  a 
counter-example to a conjecture posed in [Mi.1] and [Mi.2].

An extension of the approach from [Mi.1] was presented in a recent
paper  [MiT], in which an isomorphic version of the  stabilization
property was  investigated. Then of course it is necessary
to consider  several  variables analogues   of the above moduli,
on the expense of a clear  geometric interpretation.
This leads to a definition of upper and lower {\it envelopes\/}
(see [MiT] and also~\ref{envelopes}),
which  in  the  particular cases  of so-called bounded distortions,
give raise to the  definition  of   asymptotic $\ell_p$-spaces.
Under the same assumption of bounded distortions,
slightly different  stabilization procedure
was   also considered in [Ma].
It would be interesting to find  an isomorphic version of
stable spaces.

Notions of asymptotic type and cotype and of asymptotic
unconditionality were used in [MiSh] to study complementation.

The mentioned above  notion of an asymptotic version of a given space $X$
should   be compared with so-called  spreading model
(we recall the definition in~\ref{spreading_model}),
which also reflect some properties of a space ``at infinity''.
However, the spreading model  construction
involves only  subsequences of
a  given   sequence in $X$; thus  improving
properties of underlying space  too much, while possibly missing
phenomena occurring on  block subspaces.
In contrast,  our asymptotic
versions preserve  {all} asymptotic  finite-dimensional
properties  of a space, just choosing  its ``right''
finite-dimensional pieces
positioned  everywhere, and  then puting them together
into one infinite-dimensional  space.

\medskip

Let us  now describe  in rather imprecise terms
the intuition of an asymptotic structure of an
infinite dimensional Banach space $X$. Such a structure
is defined by a  family ${\cal B}(X)$ 
of infinite dimensional subspaces of $X$
satisfying a filtration condition which says
that for any two subspaces  from ${\cal B}(X)$ there is
a third subspace from ${\cal B}(X)$ contained in both of them
(see~\ref{filtration});
%
%
the main example  is the family ${\cal B}^0(X)$ of all
subspaces of finite codimension in $X$.  Then, for every $k$, we
define the family $\{X\}_k$ of asymptotic $k$-dimensional  spaces
associated to  this asymptotic structure  as follows
(exact  definitions are given
in the next section, see~\ref{vector_game_1}).

Fix $k$ and $\ep > 0$.  Consider a ``large enough'' number $N_1$, a
``far enough'' subspace $E_1$ of $\codim E_1 = N_1$, and an arbitrary
vector $x_1\in S(E_1)$. Next  consider a number $N_2 = N_2(x_1)$,
depending on $x_1$ and again ``large enough'',
a ``far enough'' subspace $E_2 \subset E_1$
of codimension $N_2(x_1)$ and an arbitrary vector
$x_2\in S(E_2)$.
In  the last $k$th step, we have already chosen normalized
vectors $x_1\nek x_{k-1}$ and  subspaces
$E_{k-1} \subset \cdots \subset E_2 \subset E_1$;
we then choose a ``far enough'' $E_k \subset E_{k-1}$   with
$\codim E_k = N_k(x_1\nek x_{k-1})$ and an arbitrary  vector
$x_k\in S(E_k)$. (Note that this description is
intentionally somewhat  repetitious--since
a natural meaning of
``far enough'' subspaces should imply that  their (finite) codimension
is  automatically  ``large enough''.)

We call a space $E = \spn[x_1\nek x_k] $ a permissible subspace
(up to $\ep>0$) and $\{x_i\}^k_1$---a permissible $k$-tuple
if for an {\it arbitrary\/}
choice of $N_i$ and $E_i$ (with $\codim E_i=N_i$) we would be
able to choose  normalized vectors
$\{y_i\in E_i\}$ so that a basic sequence  $\{y_i\}^k_1$ is
$(1+\ep)$-equivalent to $\{x_i\}^k_1$.

Now we  can  also clarify
the imprecise notion  of ``far enough''   subspaces  $E_i$:
by this we mean  that an {\it arbitrary\/} choice as above of
$x_i\in E_i$ results in a permissible (up to $\ep>0$) $k$-tuple
$\{x_i\}^k_1$ and a permissible (up to $\ep>0$) subspace
$E= \spn [x_i\nek x_k]$.  The existence of such subspaces ``far enough''
and of  associated $N_i$s,  will be proved in the next section
by some compactness argument.

If $F(k;\ep)$ is the set of all $k$-dimensional
$\ep$-permissible subspaces then
we put
$\{X\}_k = \bigcap_{\ep>0}F(k;\ep)$, and we call
every  space from $\{X\}_k $ a $k$-dimensional asymptotic space
of $X$.
%
%
Thus, permissible subspaces are
$(1+\ep)$-realizations of asymptotic spaces.

Finally, a Banach space $Y$ is  an asymptotic version
of $X$, if  $Y$  has  a   monotone basis
$\{y_i\}^\infty_1$ and for every $n$,  $\{y_i\}^n_1$ is a
basis in an  asymptotic space  of $X$
\ie\  $\spn[y_i]^n_1\in\{X\}_n$.

\medskip

Families of asymptotic spaces and
asymptotic versions of a given Banach  space have interesting
properties and reveal a new structure of the original space.  For
example, in Section~\ref{33} it is proved that for a fixed $p$,
with $1 \le p <\infty$, if  $X$  is a  Banach  space such that
there exists $C$  such that for every $n$, every space
$E\in\{X\}_n$ is $C$-isomorphic to $\ell_p^n$,
then every asymptotic version $Y$ of $X$ is isomorphic to
$\ell_p$ and the natural basis of $Y$ is equivalent to the
natural basis of $\ell_p$. 
It means that in such a space 
(called an {\it asymptotic $l_p$-space\/})
all permissible subspaces lie only along its
natural $l_p$ basis.

Some properties of  families of asymptotic spaces
$\{X\}_n$ can be demonstrated  through the notion of
envelopes. For any sequence  with  finite support $a \in c_{00}$
the {\it upper envelope\/} is a function
$r(a) = \sup \|\sum_ia_i e_i\|$, where
the supremum is  taken over
all natural bases  $\{e_i\}$ of asymptotic
spaces $E  \in \{X\}_n$  and  all  $n$.
Similarly, the {\it lower envelope\/} is a function
$g(a) = \inf \|\sum_ia_i e_i\|$, where
the infimum is taken over the same set.
The functions $r$ and $g$ are always very close to 
some $l_p$- (and  $l_q$)-norms (see \ref{envelopes}
for an exact statement).

An interesting general property of asymptotic versions is that some of
them are, in a sense, stable under iteration.  Precisely, we show in
Section~\ref{22} that for an arbitrary space $X$ there is a special
asymptotic version $Y$, called {\it universal}, such that its
asymptotic structure is the same as for $X$. In particular this
implies that not every space $X$, even with an unconditional basis,
can be a universal asymptotic version of any Banach space.

In Section~\ref{55} we study a complementation problem, and again the
asymptotic approach significantly simplifies the picture with respect
to ``classical" facts.


\section{Asymptotic and permissible  spaces}
\label{11} 
We  follow  [LT.1] for  standard notation  in the Banach space theory;
in particular, fundamental techniques concerning  Schauder basis,
which will be repeatedly used throughout the paper, 
can be found in    [LT.1] 1.a.

Let $X$ be a Banach space. 
By $\cB^0(X)$  we denote the family  of all
subspaces of $X$ of finite-codimension.
If  $\{u_i\}$ is  a basis  in $X$, 
or more generally,  a minimal system in $X$,
by $\cB^t(X)$
we denote the family  of all tail subspaces of $X$,
\ie\   subspaces of the form
$X^n = \overline{\spn} \{u_i\}_{i > n}$,
for some $n \in \NN$. 
 
By  ${\cal M}_n$  we denote the space of all $n$-dimensional Banach
spaces with normalized bases whose basis constant is smaller 
than or equal to 2. 
Given two such  spaces $E$, with the basis $\{e_i\}$ and
$F$, with the basis $\{f_i\}$, 
by $\di_b(E, F)$ we denote the equivalence constant
between the bases, \ie\  
$$
\di_b(E, F) = \|I: E \to F\|\, \|I^{-1}: F \to E\|,
$$
where $I$ is defined by $I e_i = f_i$, for $i=1, \ldots, n$.
Then $\log \di_b$ is a metric on   ${\cal M}_n$ which makes it into
a compact space. 

\subsection{}
\label{filtration}

An asymptotic structure of $X$ will be defined with respect
to a fixed  family $\cB(X)$ of  infinite-dimensional subspaces 
of a space $X$,  which satisfies the   filtration condition
$$
\mbox{\rm For\ every\ } X_1, X_2 \in {\cB}(X)
  \mbox{\rm \  there\  exists \ }
  X_3 \in {\cB}(X) \mbox{\rm \  such\  that\ } X_3 \subset X_1 \cap X_2.
$$
By far most important 
examples of such a family are  $\cB^0(X)$ and $\cB^t(X)$.

\subsection{}
\label{as_game}
We will work with   {\it asymptotic games\/} 
in which there are two players
\pS\   and \pV. Rules of moves are the same for all games. 
Set $X_0 = X$. In the $k$th move, player \pS\    chooses
a subspace $X_k \in \cB(X)$, and then player \pV\  
chooses a vector $x_k \in S(X_k)$ in such a way that
the vectors $x_1, \ldots, x_k$ form a basic sequence with
the basis constant smaller than or equal to 2.
Further rules will ensure that the games will stop after
a finite number of steps.

\subsection{}
\label{vector_game}
Given a space $E \in {\cal M}_n$ with a  basis $\{e_i\}$,
and  $\ep >0$, the {\it vector game\/} associated to $E$
is an asymptotic game in which the vector player \pV\  
wins  if after $n$ moves the vectors $\{x_i \}$
are $(1+\ep)$-equivalent to $\{e_i\}$.
We say that \pV\   has  a
winning strategy  for $E$ and $\ep $, if \pV\   can win every 
vector game as above.

\subsubsection{}
\label{vector_game_comments}
Since choosing by \pS\ a smaller subspace puts
\pV\ in a worst position then  the filtration property
of $\cB$ implies that without loss of generality
we can assume that, additionally, $X_{k} \subset X_{k-1}$,
for $1 \le k \le  n$.

Similarly, given $\delta >0$, by an appropriate choice of
subspaces (\cf~[LT.1] 1.a.5),  \pS\ can always ensure that 
the vectors  $\{x_i \}$ have the basis constant less than
$1 + \delta $.

\subsubsection{}
\label{vector_game_formula}

It follows that  \pV\   has  a
winning strategy for a vector game  for $E$ and
for every $\ep >0 $  if and only if
$$
\sup_{X_1} \inf_{x_1\in S(X_1)} \sup_{X_2}\ldots
\inf_{x_n\in S(X_n)} d_b([x_1,\ldots,x_n], E) = 0,
$$
with $X_k \in \cB(X)$, and $X_{k}\subset X_{k-1} $
for $1 \le k \le n$.
(Similar formulas were  used to define  fundamental notions in [MiT],
which were based on some of concepts introduced in [Mi.1].)

\subsubsection{}
\label{vector_game_1}
\begin{dfn}
A space $E \in {\cal M}_n$ with a  basis $\{e_i\}$
is called an  {\it asymptotic space\/} for $X$ if \pV\   has a winning
strategy for a vector game in $X$ for $E$ and for every $\ep >0$.
Vectors $\{x_1, \ldots, x_n\}$ in $X$
resulting from a vector game (for some asymptotic space $E$
and for  $ \ep >0 $) in which \pV\ wins,
are called a
{\it permissible} $n$-tuple and the subspace
$\spn [x_i]$ is called a  {\it permissible subspace\/}  of $X$.
\end{dfn}

So a {\it permissible subspace\/} is a $(1 + \ep)$-realization
in $X$ of  an {\it asymptotic space\/} (for some $\ep >0$).

The set of all $n$-dimensional   asymptotic spaces for $X$
is denoted by  $\{X\}_n $.
Every $E \in \{X\}_n $ has the natural basis which is monotone
(by the last comment in~\ref{vector_game_comments}).
It is easy to see that  the set $\{X\}_n $ is closed in
${\cal M}_n$.

\subsection{}
\label{subspace_game}
Given set ${\cal F}\subset {\cal M}_n$  and $\ep >0$,
the  {\it subspace game\/} is a game in which
the subspace player \pS\   wins  if after $n$ moves,
vectors $\{x_i\}$ resulting from the game
are $(1+\ep)$-equivalent
to  the basis in  some space from  ${\cal F}$.

Player \pS\   has  a winning strategy
for a subspace game  for ${\cal F}$ and $\ep $,
if \pS\   can win  every  such game.
The filtration property clearly implies
that \pS\ can always choose subspaces satisfying
$X_{k} \subset X_{k-1}$ for $ 1 \le k \le n$,
this way only improving
his chances to win. Therefore we will always assume
that  winning strategy for \pS\ in a subspace game
satisfies this condition.

\subsubsection{}
\label{subspace_game_formula}
It follows that player \pS\   has  a
winning strategy for a subspace game  for $\cal F$ and
for every $\ep >0 $  if and only if
$$
\inf_{X_1} \sup_{x_1\in S(X_1)} \inf_{X_2}\ldots
\sup_{x_n\in S(X_n)} \inf_{F \in {\cal F}}
 d_b([x_1,\ldots,x_n], F) = 0,
$$
with $X_k \in \cB(X)$, and $X_{k}\subset X_{k-1} $
for $1 \le k \le n$.

\subsubsection{}
\label{subspace_game_tree}
To  visualize both  formulae~\ref{vector_game_formula}
and~\ref{subspace_game_formula},
we can think about
a certain tree-like structure of subspaces $X_k$
from the partially ordered by inverse inclusion  set $\cB(X)$,
and arbitrary vectors  $x_k \in S(X_k)$, with choices of
subsequent subspaces depending on the earlier vectors.
Then the vector player \pV\ has a winning strategy in
a vector game for  some  space $E \in {\cal M}_n$ and $\ep >0$,
if \pV\  can find, {\it arbitrarily far\/} along  $\cB(X)$,
vectors $\{x_i\}$ which are $(1+\ep)$-equivalent
to the basis in $E$.
The subspace player \pS\ has a winning strategy
in a subspace game for  a  subset $\cal F \subset {\cal M}_n$
and  $\ep >0$,
if by choosing subspaces $X_k$ {\it far enough\/} along
$\cB(X)$, \pS\  can ensure  that the  vectors $\{x_k\}$
are  $(1+\ep)$-equivalent  to the basis of some
space from  $\cal F$.

\subsubsection{}
\label{subspace_game_warning}
We will show in~\ref{games} below that
the subspace player
\pS\   has a winning strategy in a subspace game for
$\{ X\}_n$ and for every $\ep >0$.
We will then repeatedly  use this fact to show
that  if in an arbitrary asymptotic game
\pS\ follows his winning strategy for a fixed $\ep >0$,
then a subspace of $X$ resulting in the game
is a $(1 + \ep)$-representation of some
asymptotic space from $\{X\}_n$.
This yields in particular
that this subspace is permissible,
without actually  stating  which space from $\{X\}_n$
does it represent.

\subsubsection{}
\label{subspace_game_1}
Let $\cal W$ be the family of all  closed subsets
${\cal F}$ of ${\cal M}_n$ such that \pS\   has a winning strategy
in a subspace game for ${\cal F}$ and for every $\ep >0$.
Clearly, ${\cal M}_n \in {\cal W}$ and
$\emptyset \not\in {\cal W}$.
Moreover, the filtration property
immediately implies
that if $ {\cal F}_i \in {\cal W}$,
for $i=1, \ldots, m$, then
$\bigcap_i {\cal F}_i \in \cal W$.
Let $\widetilde{\cal F}
= \bigcap_{{\cal F} \in {\cal W}} {\cal F}$.
This is a non-empty closed  subset of ${\cal M}_n$.
We shall show that $\widetilde{\cal F} \in {\cal W}$.

This follows from a compactness argument.
Let ${\cal D}_\delta (F)$ be the open ball in ${\cal M}_n$ of radius
$\delta$ and center at $F$.  Observe that for $\delta >0$, the set
${\cal F}^\delta
= \bigcup_{F \in \widetilde{\cal F}} {\cal D}_\delta (F)$
contains an intersection of a finite number of sets from $\cal W$.
Indeed, the complement $({\cal F}^\delta)^c$
of ${\cal F}^\delta$ is compact, and
it is contained in $ (\widetilde{\cal F})^c$, which in turn
is covered by the union of complements of sets from $\cal W$.
Thus for every  $\delta >0$,
${\cal F}^\delta$ contains a set from $\cal W$,
hence  \pS\   has a winning strategy
in a subspace game for ${\cal F}^\delta$  and every $\ep >0$.
Since these sets approximate
$\widetilde{\cal F}$ arbitrarily close,
for an arbitrary fixed $\ep >0$,
\pS\   has a winning strategy
in a subspace game for $\widetilde{\cal F}$
as well. Thus  $\widetilde{\cal F} \in {\cal W}$.

\subsection{}
\label{games}
The set of asymptotic spaces  $\{ X\}_n$ coincides with
$\widetilde{\cal F}$.
Therefore the subspace player
\pS\   has a winning strategy in a subspace game for
$\{ X\}_n$ and for every $\ep >0$.
In particular, $\{ X\}_n$ is non-empty.

First,  $\{ X\}_n \subset \widetilde{\cal F}$.
Indeed, let $E \in \{X\}_n$ and let $\ep >0$.
Consider an asymptotic game
in which each player follows his own
strategy;
player \pS\   follows the winning strategy for
a subspace game for
$\widetilde{\cal F}$, and
player \pV\   follows the winning strategy for
a vector game for $E$.
Strategy of \pV\   implies that vectors $\{x_i\}$
resulting from this game are $(1+ \ep)$-equivalent to the
basis in $E$; strategy of \pS\   implies that
they are also $(1+ \ep)$-equivalent to the
basis in some space from $\widetilde{\cal F}$.
Hence $\di_b (E, \widetilde{\cal F}) \le (1 + \ep)^2$,
for every $\ep >0$. Thus $E \in \widetilde{\cal F}$,
since  $\widetilde{\cal F}$ is closed.

Next, observe that if $E \not \in \{X\}_n$,
then  for  $\ep_0 >0$ sufficiently small,
player \pS\   has a  strategy in a subspace
game such that resulting vectors $\{x_i\}$
satisfy
$\di_b (\spn[x_i], E)\ge 1 + \ep_0$.
Thus for every $\ep < \ep_0/2$,
\pS\   has a winning strategy in a subspace game
for  ${\cal F}' = \widetilde{\cal F}
\backslash  {\cal D}_{\ep_0/2} (E) $.
Since ${\cal F}'$ is closed, the minimality
of  $\widetilde{\cal F} $ implies
in particular that
$E \not\in  \widetilde{\cal F}$.


\subsection{}
\label{remarks_permiss_1}
We look in more detail at the  family  
of all asymptotic spaces of $X$.

Spaces $l_p$ play a special role here.
For $1 \le p < \infty$, the standard unit vector basis
in $l_p$ is denoted by $\{\mbox{\bf e}_i\}$. The same notation is
used in $c_0$ and in finite-dimensional spaces $l_p^n$.

\subsubsection{}
\label{compar_families}
Denote by $\{X\}^0_n$ and by $\{X\}^t_n$  the sets of all
$n$-dimensional asymptotic spaces with respect to the families
$\cB^0(X)$ and $\cB^t(X)$.
Clearly,  $\{X\}^0_n \subset \{X\}^t_n$.
In general, this inclusion is proper; however
the main property of  shrinking systems  immediately   implies
that  if a fundamental system in $X$  is  shrinking,
then  $\{X\}^0_n = \{X\}^t_n$.

As an example of a space for which asymptotic structures
depend on  a  family $\cB(X)$, consider
the space $c$ of all convergent scalar sequences.
Let $\{u_i\}$ be the natural basis  in $c$, that is,
$u_1 = (1, 1, 1, \ldots)$
and $u_i = \mbox{\bf e}_{i-1}$ for $i >1$,
and  consider $n$-dimensional asymptotic  spaces with respect to
the family of  tail  subspaces of   $\{u_i\}$.
It is obvious that the only such
space  is $l_\infty^n$ with the  standard unit vector basis.
This in particular implies that
$\{c\}^0_n = \{l_\infty^n\}$.
On the other hand, consider a conditional basis
$\{v_i\}$ in $c$ given by
$v_i = \sum_{j=i}^\infty \mbox{\bf e}_{j}$
for $i=1, 2, \ldots$ and consider the set
$\{c\}^t_n$ with respect to this basis.
Clearly, $l_\infty^n \in \{c\}^t_n$,
however  it is  easy to see that this set
is larger: it also contains
the space $E $ which is $l_\infty^n$ with the {\it conditional\/}
basis  $\widetilde{v}_i = \sum_{j=i}^n \mbox{\bf e}_{j}$
for $i=1, \ldots, n$.

\subsubsection{}
\label{spreading_model}
Recall that  a bounded  non-convergent sequence  $\{z_i\}$ in a
Banach space $X$ is said to generate spreading model,  if for every
finite sequence of scalars  $(a_1, \ldots, a_k)$
the $k$-fold limit
$ \lim_{n_1}\ldots \lim_{n_k} \|\sum_{i=1}^k a_i z_{n_i}\|$
exists, as $n_i\to \infty$, for $i=1, \ldots, k$,
with $n_1 < \ldots < n_k$.

Then one can define  the spreading model $F$  as a Banach space
with the basis $\{\mbox{\bf f}_i\}$
such that for every finite sequence of scalars
$(a_i)$  one has
$$
\|\sum_{i=1}^k a_i \mbox{\bf f}_i\| = \lim_{n_1}\ldots \lim_{n_k}
       \|\sum_{i=1}^k a_i z_{n_i}\|.
$$
Clearly,  the basis   $\{\mbox{\bf f}_i\}$  is
spreading invariant, \ie\ for every finite sequence of scalars
$(a_i)$  and every $n_1 < n_2 < \ldots $ one has
$\|\sum_{i} a_i \mbox{\bf f}_i\| =
\|\sum_{i} a_i \mbox{\bf f}_{n_i}\|$.
In such a situation,
the sequence of differences
$\{\mbox{\bf f}_{2i} - \mbox{\bf f}_{2i-1} \}$ is unconditional
(and clearly still spreading invariant).

It is a well-known result by Brunel and Sucheston [BS]
and it follows from Ramsey's theorem that
every bounded sequence with no Cauchy subsequences
contains a subsequence generating   spreading model;
and then  the differences of this subsequence generate
an unconditional spreading model
(\cf~also [MiSch] Section 11).
The reader can consult \eg~[BL] on more details on spreading models.

\subsubsection{}
\label{krivine}
Recall that if a sequence  $\{z_i\}$ generates  an
unconditional spread\-ing mod\-el, then a direct application
of Krivine's theorem [K] says that there exists $1 \le p \le \infty$
such that for every $n \in \NN$ and every
$\ep >0 $ there is a finite  scalar sequence
$\al=\{\al_1, \ldots, \al_m\}$ such that
any $n$ successive blocks $\{x_j\}$ of  $\{z_i\}$
with the same distribution $\al$ and ``far enough'',
are $(1+\ep)$-equivalent to the basis $\{\mbox{\bf e}_i\}$
in $l_p^n$.

Since every Banach space has a sequence
generating unconditional spreading model,
there exists $ 1 \le p \le \infty$ such that
$l_p^n \in \{X\}_n$ for every $n$.

\subsubsection{}
\label{stabilized}
Let us briefly discuss a concept of a stabilized asymptotic
structure, which appears implicitely or explicitely
in many papers already mentioned
([G.1], [Ma], [MiT]) and others ([C], [G.2]).
This concept allows for passing to infinite-di\-men\-sion\-al
subspaces of a given space $X$ and hence we can assume that
$X$ has a basis.

By $\cB_\infty(X)$ denote  the set  of all
infinite-dimensional block subspaces of $X$.
We can then consider a family ${\cal D} =
\{{\cal D}_n\}$ of subsets  ${\cal D}_n \subset {\cal M}_n$,
for $n=1, 2, \ldots$, such that
there exists $Y \in \cB_\infty(X)$ such that for every
$n \in \NN$   the following stabilization condition holds:
for   every $Z \in \cB_\infty(Y)$,
we have  $ {\cal D}_n = \{Z\}_n$.
It is not difficult to show that the sets
${\cal D}_n $ are non-empty.
In fact,  there exists $ 1 \le p \le \infty$ such that
$l_p^n \in {\cal D}_n$ for every $n$.
Each space $E$ from  ${\cal D}_n $
(for $n \in \NN$)  is called a
{\it stabilized  asymptotic space\/}  for $X$,
and $Y$ is called a stabilizing subspace for ${\cal D} $.

To show that all the ${\cal D}_n $'s are non-empty and to construct $Y$,
first observe that the compactness of ${\cal M}_n$ and Zorn's lemma
show that  for a fixed $n$, given $\widetilde{Y} \in \cB_\infty$,
there exists $Y_n \in \cB_\infty( \widetilde{Y} )$ such that
the $n$th stabilization condition holds in $Y_n$. Then the  space
$Y$ is a diagonal  subspace of $Y_n$'s.
The last statement about  $l_p^n$'s follows from the fact
that the set of  all $p$'s such that
$l_p^n \in \{X\}_n$ for every $n$ (\cf~\ref{krivine}),
is closed.

\subsubsection{}
\label{gowers}
An important  recent combinatorial theorem by Gowers [G.1]
provides  further general information on  families
of stabilized asymptotic spaces.
Let $X$ be a Banach space with a basis.
Let  $\Sigma$ be a set of all  sequences
$\{x_1, \ldots, x_n\}$, where $n \in \NN$
and the vectors are successive normalized  blocks
of the basis.
A subset  $\sigma \subset  \Sigma$ is called
{\it large\/} if for every $Y \in \cB_\infty(X)$,
there is a sequence  $\{x_1, \ldots, x_n\}\in \sigma$
with $x_i \in Y$ for  $i =1, \ldots, n$.

Given  a subspace  $Y \in \cB_\infty(X)$,
consider a {\it general infinite-dimensional vector game\/}
inside  $Y$,
which is in essential way  less restrictive than the game
introduced in~\ref{vector_game}.
Here for  the $k$th move of the game,  the subspace player
\pS\ chooses  a  subspace $Y_k \in \cB_\infty(Y)$
and then the vector player \pV\  chooses
a vector $x_k \in S(Y_k)$.
Given a set  $\sigma \subset  \Sigma$ and  $\ep >0$,
player \pV\ wins the game inside $Y$, if  after some
number of moves
the sequence $\{x_1, \ldots, x_n\}$ he has chosen
is  $(1+\ep)$-equivalent to a sequence from $\sigma$.
Note, in comparison with~\ref{vector_game} and~\ref{as_game},
that here subspaces chosen by \pS\
may have {\it infinite\/} codimension
and the  number of moves in the game is not prescribed
in advance; in fact, this number even does not have to be finite,
if \pV\  does not have a winning strategy.

Gowers' theorem says that if $\sigma \subset  \Sigma$ is large
on $X$ then for every $\ep >0$ there is a subspace
$Y \in \cB_\infty(X)$ inside which \pV\ has a winning strategy
for $\sigma$ and $\ep$.

Now let  $X$  be an arbitrary space and let $n \in \NN$.
If  $E \in {\cal M}_n$
and there  is  an infinite-di\-men\-sion\-al  subspace
$Z \subset X$ with a basis
such that for every  $Y \in \cB_\infty(Z)$
and every $\ep >0$, there are $n$
successive blocks $\{v_1, \ldots, v_n\}$  in $Y$
such that $\di_b(\spn [v_i], E) \le 1 + \ep$,
then $E $ is a stabilized  asymptotic space for
$X$  and in particular,
$E\in \{X\}_n$.

Indeed, given $\ep >0$,  we use Gowers' theorem
for the set  $\sigma$ of all $n$-tuples of successive blocks
$\{v_1, \ldots, v_n\}$ as above; this  $\sigma$
is large  on $Z$.

Notice that this argument does not require the full strength of
Gowers' result: the game used above has a fixed length and in this
case the theorem is easier.
Let us also mention that as an easy corollary to his general result,
Gowers obtained the following attractive structure dichotomy for
Banach spaces: every infinite-dimensional Banach space  either has
a subspace with an unconditional basis or has a hereditarily
indecomposable subspace $X_0$ (\ie\ no subspace of $X_0$ is a topological
direct sum $Y \oplus Z$ of infinite-dimensional subspaces.)

\subsubsection{}
\label{odell-schl}
Recall a recent  construction of Odell and Schlumprecht
[OS.1], of a Banach space $Z$ with a basis $\{z_i\}$ such that
for every $n \in \NN$,  every $n$-dimensional space
$E $ with a monotone basis,
every  $Y \in \cB_\infty(Z)$
and every $\ep >0$, there are $n$
successive blocks $\{v_1, \ldots, v_n\}$  in $Y$
such that $\di_b(\spn [v_i], E) \le 1 + \ep$.
It follows from~\ref{gowers} that for this space $Z$,
every finite-dimensional space  with a monotone basis
is a stabilized asymptotic space.

\subsection{}
\label{asymptotic_l_p}
\begin{dfn}
A Banach space   $X$ is called an {\it asymptotic}-$l_p$ space,
for $ 1 \le p \le \infty$,
if there is a constant $C$ such that for every $n$ and every
$E \in \{X\}_n$ we have $\di_b (E, l_p^n) \le C$.
The asymptotic structure $\{X\}_n$ is
determined by the family $\cB^0(c)$ of  all finite-codimensional
spaces of $X$, \ie\  $\{X\}_n= \{X\}_n^0$.

\end{dfn}

An example discussed in~\ref{compar_families} shows
that the restriction of the asymptotic
structure to $\{X\}_n^0$  is essential in general:
the space $c$ is an asymptotic-$l_\infty$, but
{\it some\/} asymptotic spaces
relative to the family $\cB^t(c)$ of tail spaces of the
conditional asis
$\{v_i\}$, are not equivalent (in sense of $\di_b(\cdot, \cdot)$)
to $l_\infty^n$ with the standard basis.

\subsubsection{}
\label{asymptotic_l_p_1}
If $l_p^2$ is the only  2-dimensional asymptotic space for
a Banach space $X$, \ie\ $\di_b (E, l_p^2) =1$,
for every $E \in \{X\}_2$, then
$X$ contains almost isometric copies of $l_p$.
Indeed,  a  well-known easy  argument  shows that
the formula from~\ref{subspace_game_formula}
allows to construct, for every $\ep >0$,
a basic sequence $\{x_i\}$ in $X$  such that
$\{x_i\} \stackrel{1+\ep}\sim \{\mbox{\bf e}_i\}$
(\cf~\eg~[MiT]).

\subsubsection{}
\label{asymptotic_l_p_2}
A trivial example of  an asymptotic-$l_p$   space
not isomorphic to $l_p$,
is an $l_p$-direct sum  $(\sum \bigoplus Z_l)_p$
of finite-dimensional spaces.

A class of much more sophisticated examples  are
$p$-convexified Tsirelson spaces $T_{(p)}$;
these spaces are   asymptotic-$l_p$ and  they do not
contain subspaces isomorphic to $l_p$ (\cf~\eg~[CS]).

\subsubsection{}
\label{asymptotic_l_p_stab}
If  all $n$-dimensional  stabilized
asymptotic spaces (\cf~\ref{stabilized})
are uniformly equivalent to the unit
vector basis in  $l_p^n$, for $n \in \NN$,
and $X$ itself is a stabilizing subspace,
then $X$ is called a  {\it stabilized  asymptotic}-$l_p$
space.  These spaces were  investigated
in [MiT] and [Ma] (where they were called
just asymptotic-$l_p$ spaces).

\subsubsection{}
\label{asymptotic_l_p_3}
From the point of view of Banach space theory
it is tempting to consider a seemingly more general
concept than  asymptotic-$l_p$  spaces, in which
the condition that
the basis  in $E$ is $C$-equivalent to the natural basis
in $l_p^n$, is replaced by the condition that
$E$ itself is  $C$-isomorphic  to $l_p^n$
(for $E \in \{X\}_n$).
Recall that a Banach space with a basis has
uncountably many mutually non-equivalent
bases (\cf~[LT.1] 1.a.8); for spaces $l_p$, with $ 1 < p < \infty$,
$p \ne 2$,
these bases may be chosen to be  even unconditional
(\cf~[LT.1] 2.b.10).
It is therefore rather striking that
in the asymptotic setting  discussed here
for $ 1 \le p < \infty$,
the more general  condition  of isomorphism
of asymptotic spaces to $l_p^n$
already implies  the equivalence  of the natural bases.
This will be  proved in  Section~\ref{33}.

\subsection{}
\label{remarks_permiss_2}
Let us  consider again an asymptotic structure with respect to an
arbitrary family $\cB$ satisfying  the  filtration
condition~\ref{filtration}.
We shall discuss some properties of asymptotic
families which show an interplay between  different level
families.

\subsubsection{}
\label{isomorphic}
Let $X^{(1)}$ and $X^{(2)}$ be two $C$-isomorphic Banach spaces.
For every $n \in \NN$, the Hausdorff distance (in ${\cal M}_n$) between
$\{X^{(1)}\}_n$  and  $\{X^{(2)}\}_n$  is smaller than or equal to $C$.
That is, if $1 \le i\ne j \le 2$ then for every
$E \in \{X^{(i)}\}_n$  there is $F \in \{X^{(j)}\}_n$
such that  $\di_b (E, F ) \le C$.
In particular, if $X^{(1)}$ is an asymptotic-$l_p$ space
then so is  $X^{(2)}$.

Indeed, let $T: X^{(1)} \to X^{(2)}$ be  an isomorphism.
Given  $E \in \{X^{(i)}\}_n$,  the corresponding space $F$
will be spanned by an $n$-tuple  resulting from
a subspace game in $ \{X^{(j)}\}_n$; together with this game
one considers a vector game  for $E$  in $\{X^{(i)}\}_n$,
and the moves between the two games are  translated one to another
by the operators $T$ and $T^{-1}$.

\subsubsection{}
\label{mixtures}
Let $n_1, \ldots, n_k$ be natural numbers.
Let  $E_j \in \{X\}_{n_j}$, for $j=1, \ldots, k$.
For every $N \ge \sum_j n_j$ and any disjoint subsets
$I_j$ of $\{1, \ldots, N\}$,
with $|I_j|= n_j$  for $j=1, \ldots, k$,
there exists an asymptotic space $F \in \{X\}_N $
with a basis $\{f_i\}$ such that
$\di_b (\spn [f_i]_{i \in {I_j}}, E_j) =1$
for $j=1, \ldots, k$.

Indeed, let $\ep >0$.
Consider an asymptotic game which ends  after $N$ moves.
Player \pS\   simply follows his winning strategy for
a subspace game for $\{X\}_N$ and $\ep$. Strategy for
player \pV\   is more  complicated.
For $j=1, \ldots, k$, write $I_j =
\{i_1^{(j)}, \ldots,  i_{n_j}^{(j)} \}$;
if $i \in I_j$ for some $1 \le j \le k$,
say $i = i_{l}^{(j)}$ for  $1 \le l \le n_j$, then
\pV\   makes his choice of vector $x_i$ following
the winning strategy for the $l$th move in a vector game
for $E_j$ and  $\ep$, as if his  previous
choices in this game were the vectors $x_s$,
for $s = i_1^{(j)}, \ldots,  i_{l-1}^{(j)}$.
If $i \not \in I_j$ for any $j$, \pV\   picks the vector
$x_i$ arbitrarily.

Consider  the   vectors $\{x_i\}$ resulting in the game.
The strategy of \pS\ implies that  they are
$(1+\ep)$-equivalent to the  basis  $\{f_i^{\ep}\}$ in some
asymptotic $N$-dimensional space
$F_\ep \in \{X\}_N$.
The strategy of \pV\ implies in turn that
$\di_b (\spn [f_i^{\ep}]_{i \in {I_j}}, E_j) \le 1 + \ep$,
for  $1 \le j \le k$.
Then a required space  $F \in {X}_N$  is any cluster point
in ${\cal M}_N$ of the $F_\ep$'s, as $\ep \to 0$.

\subsubsection{}
\label{as_blocking}
Let  $E \in \{X\}_n$, and let  $F$  be a
block subspace of $E$, that is, $F$ is spanned
by successive blocks of the basis in $E$.
Then $F$ is  an asymptotic space,
$F \in \{X\}_m$, where $m = \dim F$.
Moreover, given $n \in \NN$ and $\ep >0$, the subspace player
\pS\ has a strategy in an asymptotic game such that
after $n$ moves,  all  normalized successive blocks
of the $n$-tuple  resulting from the game,
are permissible, \ie\ each of them is $(1+\ep)$-
equivalent to the basis in some asymptotic space.

To prove the first statement,
let $\{e_i\}$ be the basis in $E$ and let $\{u_k\}$
be  successive  blocks of $\{e_i\}$ spanning $F$.
Let $i_0 = 1 < i_1 < \ldots < i_m = n+1$
such that $u_k = \sum _{i=i_{k-1}}^{i_k -1} a_i e_i$, for
$k=1, \ldots, m$.
Given  $\ep >0$,  consider   a vector game  for $E$
and $\ep$ in which  choices of  player \pS\
follow the pattern
$X_1, X_1, \ldots, X_1, X_2, \ldots, X_2, X_3, \ldots$,
with  the change of a subspace  being made only in
the $i_k$th moves and subspaces $X_k$
being arbitrary ($k=0, \ldots, m-1$),
and  player \pV\   follows his winning
strategy for $E$.
Denote  the resulting  permissible $n$-tuple by
$\{x_i\}$, then
$\{x_i\} \stackrel{1 + \ep}\sim \{e_i\}$.
Moreover, the blocks
$v_k = \sum _{i=i_{k-1}}^{i_k -1} a_i x_i$
obviously satisfy $\{v_k\} \stackrel{1 + \ep}\sim \{u_k\}$.
This describes a winning strategy for \pV\   in a vector
game for $F$ and $\ep$.  Hence $F \in \{X\}_m$.

For the moreover part, it is not difficult to see
from a simple perturbation argument, that
if \pS\ follows his winning strategy for $\{X\}_n$
and $\delta >0 $, then arbitrary successive normalized
blocks  $\{w_k\}$  of any $n$-tuple $\{x_i\}$
resulting in the game,
are $(1+\delta)(1+n\delta)$-equivalent
to  corresponding normalized blocks of the basis
in the  space  from $\{X\}_n$ associated to
$\{x_i\}$. Thus $\{w_k\}$  are permissible.

\subsection{}
\label{envelopes}
We  conclude this section by  introducing the notion of 
{\it envelopes\/}  which  is of inde\-pen\-dent interest.

\subsubsection{}
\label{env_1}
Recall that $c_{00}$  denotes the space of all scalar
sequences eventually zero.
The {\it upper\/} and the {\it lower envelopes\/} 
for $X$ are functions
$r(\cdot)$ and $g(\cdot)$, respectively, defined  for
$a = (a_1, \ldots, a_n, 0 \ldots) \in c_{00}$  by
$r(a) = \sup \|\sum_ia_i e_i\|$
and   $g(a) = \inf \|\sum_i a_i e_i\|$, where
the supremum and the infimum are taken over
all natural bases  $\{e_i\}$ of asymptotic
spaces $E  \in \{X\}_n$  and  all  $n$.

The functions  $r(\cdot)$ and $g(\cdot)$ are obviously
unconditional  and subsymmetric. It is easy to see that
$r(\cdot)$ is a norm on $c_{00}$
and that $g(\cdot)$ satisfies triangle inequality
on disjointly supported vectors.
These functions were used in an essential way in [MiT].

\subsubsection{}
\label{env_2}
Note that the upper envelope is 
{/it sub-homogeneous\/}. By this we mean that
for any finite number of successive vectors
$b^i \in c_{00}$   such that $r (b^i)\le 1$ for $i=1, 2, \ldots$
and for any vector $a = (a_i)_i\in c_{00}$,
we have
$$
r(\sum_i a_i b^i ) \le r(a).
$$
Similarly, the lower envelope satisfies 
the {/it super-homogeneity\/} condition:
if $ g(b^i) = 1$ for $i=1, 2, \ldots$ then
$$
g(\sum_i a_i b^i ) \ge g(a).
$$
The proof of both inequalities usess~\ref{as_blocking}
and unconditionality of both functions.

\subsubsection{}
\label{env_3}
It is a general and  interesting fact that 
sub-homogeneous norms or functions satisfying a weaker
triangle inequality as  $g(\cdot)$ does, are always
close to some $l_p$- norm. We formulate the exact statement
for our envelope functions.

There exist $1 \le p, q \le \infty$  and $C, c >0$ and
for every $\ep >0$ there exist $C_\ep, c_\ep >0$ such  that
for $a \in c_{00}$ we have
$$
 c_\ep \|a\|_{l_{q+\ep}} \le g(a) \le C \|a\|_{l_{q}}
          \quad \mbox{and} \quad
 c\|a\|_{l_p} \le r(a) \le C_\ep \|a\|_{l_{p-\ep}}.
$$

We outline a  standard argument for the function $r(\cdot)$.
%
%
For a   positive integer $n$ set
$\la_r(n) = r((1, \ldots, 1, 0,\ldots))$.
Then  sub-homogeneity  of $r(\cdot)$ 
discussed in~\ref{env_2} implies that
$\la_r (n\, m) \le \la_r (n)\, \la_r ( m)$. By induction,
we get $\la_r (n^k) \le \la_r (n)^k $.
Let $ 1/p = \inf \ln \la_r(n)/ \ln n$. Clearly,
$\la_r (n) \ge n^{1/p}$ for all $n$. On the other hand,
for every $\ep >0$ there exists a constant $C_\ep$ such that
$\la_r(n) \le C_\ep n^{1/(p - \ep)}$.

By Krivine's theorem for the space $(c_{00}, r(\cdot))$, this space
contains $l_p^n$'s uniformly on succesive blocks of the natural
basis. Using  submultiplicativity of $r$, it easily follows
that $r(a) \ge c\|a\|_{l_p}$, for all $a \in c_{00}$.
On the other hand, it can be also easily seen that
an upper power type estimate for $\la_r(n)$ implies
a similar estimate for $r$ (with different $C_\ep$, though).


\section{Asymptotic versions}
\label{22}
In this section we 
introduce   infinite-dimensional
spaces which reflect properties of the whole 
sequence  $\{\{X\}_n\}$ of families  of $n$-dimensional asymptotic 
spaces of  a given Banach space $X$.
This will be done by considering an additional
structure given naturally by an inclusion
on bases of asymptotic  spaces.

\subsection{}
\label{as_version}
A space $Y$ with a  monotone basis $\{y_i\}$
is called an {\em asymptotic version\/} of $X$
if for every $n \in \NN$ we have
$\{y_i\}_{i=1}^n  \in \{X\}_n$.
The set of all asymptotic versions of $X$ is denoted by
${\cal A} (X)$.

A construction of  an asymptotic version of a given
space $X$, fully resembles  the concept of 
an injective  limit.
First observe   that if 
$\{f_i\}_{i=1}^n   \in \{X\}_n$ 
then the restriction $\{f_i\}_{i=1}^{n-1}$
is in  $\{X\}_{n-1}$.

Conversely, for every 
$\{e_i\}_{i=1}^n \in    \{X\}_n$ 
there is 
$\{f_i\}_{i=1}^{n+1}    \in  \{X\}_{n+1}$ 
such that 
$\{f_i\}_{i=1}^n \stackrel{1}\sim \{e_i\}_{i=1}^n $.

Indeed, given $\ep >0$, consider an asymptotic game 
in $X$ which ends after
$n+1$ moves; in which player \pS\    follows his winning startegy
in a subspace game for  $\{X\}_{n+1}$, 
and player \pV, in the first $n$ moves,
follows a winning strategy  in a vector game
for  $\{e_i\}_{i=1}^n $, and in the $(n+1)$th move 
picks an arbitrary vector.  Denote the  resulting
$(n+1)$-tuple by $\{f_i^{\ep}\}_{i=1}^{n+1}$. 
An argument similar to the one used at the end
of~\ref{mixtures} shows that  any
cluster point of the  $\{f_i^{\ep}\}_{i=1}^{n+1}$'s in 
${\cal M}_{n+1}$  belongs to $\{X\}_{n+1}$ and 
its restriction is clearly  $\{e_i\}_{i=1}^{n} $.

We can  construct an increasing sequence
$F_1 \subset \ldots \subset F_n  \subset F_{n+1} \subset \ldots$
with bases
$\{f_1\} \subset \ldots \subset \{f_i\}_{i=1}^n \subset
\{f_i\}_{i=1}^{n+1}  \subset \ldots $
such that $F_n \in \{X\}_n$.
Then $Y = \overline{ \bigcup_n F_n} $ 
is an asymptotic version of $X$.

\subsection{}
\label{example_as_version}
Let us consider few simple examples of asymptotic versions.

Clearly, a space $X$ is an asymptotic-$l_p$  if and only
if all asymptotic versions of $X$  are uniformly equivalent
to the standard unit vector basis in $l_p$.

The space $Z$ from~\ref{odell-schl} has every Banach space 
$Y$ with a monotone basis as its asymptotic version.

\subsubsection{}
\label{uncond_as_version}
A space  $X$ has an {\it asymptotic unconditional  structure\/}
if   there exists $C$ such that for every  
asymptotic space $E  \in \{X\}_n $ (where $n = \dim E$)
the natural basis $\{e_i\}$ in $E$ is $C$-unconditional,
\ie\  $\mbox{unc\,}\{x_i\}\le C$.

Clearly, $X$ has an  asymptotic unconditional  structure
if and only if  there exists $C$ such that for
every asymptotic version $Y$ of $X$ 
the natural basis in $Y$  is $C$-unconditional.

\subsection{}
\label{univ_as_version}
Let  $Y \in {\cal A}(X)$, let $n \in \NN$ and 
let  $E \in \{Y\}_n$. Fix $ \ep >0$. Then the basis
in $E$ is $(1 +\ep)$-equivalent to $n$ successive blocks
of  some initial interval of the basis in $Y$, say
$\{y_i\}_{i=1}^N$. Since $\{y_i\}_{i=1}^N \in \{X\}_N$,
then~\ref{as_blocking} implies that $E$ is
$(1 +\ep)$-close to some asymptotic space  for $X$.
Thus $\{Y\}_n \subset \{X\}_n$,   for every $n \in \NN$.

The following theorem shows that we can construct
an asymptotic version of $X$ which contains  all
asymptotic spaces of $X$ in an asymptotic way.

\begin{thm}
For every Banach space $X$ there  exists an asymptotic version
$Y \in {\cal A}(X)$  such that $\{Y\}_n = \{X\}_n$
for every $n \in \NN$. Moreover, $Y$ can be constructed
in such a way that every asymptotic space of $X$ is represented
(in an asymptotic way) as a permissible span of basic vectors of $Y$.
\end{thm}

Such a space $Y$ is  called a {\em universal asymptotic version\/}
for $X$.

\subsection{}
\label{univ_corol}
It follows that  not every Banach space
can be a universal asymptotic version of another Banach space.
Examples from~\ref{example_as_version} imply that this is
a case of  an asymptotic-$l_p$ space not
isomorphic to $l_p$ (see~\ref{asymptotic_l_p_2}),
or of a space with an asymptotic unconditional
basis which is not  unconditional ([G.2]).

\subsection{}
\label{proof}
The proof of Theorem~\ref{univ_as_version} is based on several lemmas.

\subsubsection{}
\label{lemma_1}
The first lemma is similar to~\ref{mixtures} and has
an analogous proof which is left for the reader.
\begin{lemma}
Let $n, m \in \NN$, let $N \ge n\,m$.  Let
${\cal I} = \{I_j \}_{j=1}^m$
be a family  of $m$ subsets of  $\{1, \ldots, N\}$,
such that $|I_j|= n_j\le  n$  for $j=1, \ldots, m$,
and the following condition is satisfied:
for arbitrary two sets $I_k$ and $I_l$ in $\cal I$,
the intersection $ I_k \cap I_l$
is either empty or it is an initial interval
of each of them, \ie\   if
$I_k = \{t_1, \ldots, t_{m_k}\}$
and   $I_l = \{s_1, \ldots, s_{m_l}\}$,
and if  $t_\mu = s_\nu$ for some $\mu, \nu \in \NN$,
then  $\mu = \nu$ and $t_1 = s_1$, $\ldots$, $t_\mu = s_\mu$.
Let $E \in \{X\}_n$ with a basis $\{e_l\}$. There exists
an asymptotic space $F \in \{X\}_N $
with a basis $\{f_i\}$ such that
$$
\{f_i\}_{i \in {I_j}} \stackrel{1}\sim \{e_l\}_{l=1}^{n_j},
$$
for $j=1, \ldots, m$.
\end{lemma}

\subsubsection{}
\label{lemma_2}
We also require infinite-dimensional facts of a similar
nature.  To avoid unnecessary repetitions, let us
use the convention that if a basis  $\{z_i\}$
of a Banach space $Z$ is understood from the context,
for a  basic  sequence $\{y_i\}$ we shall write
$\{y_i\} \stackrel{1} \sim Z$ instead of
$\{y_i\} \stackrel{1} \sim \{z_i\}$.

The proof of the next lemma follows by combining~\ref{mixtures}
and~\ref{as_version}.

\begin{lemma}
Let $Y_1$ and $Y_2$ be two asymptotic versions of $X$.
Let $I_1$ and $I_2$ be two infinite disjoint subsets of $\NN$.
There exists an asymptotic version $Y$ of $X$ with a basis $\{y_i\}$
such that   $\{y_i\}_{i \in {I_1}} \stackrel{1}\sim Y_1 $
and   $\{y_i\}_{i \in {I_2}} \stackrel{1}\sim Y_2$.
\end{lemma}

\subsubsection{}
\label{lemma_3}
The final lemma is a version of the latter one for infinitely
many spaces. We leave the proof to the reader.

\begin{lemma}
Let $\{Y_j\}$ be a sequence of asymptotic versions of $X$.
Let $\{I_j\}$  be a sequence of infinite mutually disjoint
subsets of $\NN$.
There exists an asymptotic version $Y$ of $X$ with a basis $\{y_i\}$
such that   $\{y_i\}_{i \in {I_j}} \stackrel{1}\sim Y_j $
for every $j = 1, 2, \ldots$.
\end{lemma}

\subsubsection{}
\label{real_proof}
Now we are ready for the proof of the theorem.

\proof
Fix an arbitrary asymptotic space $E \in \{X\}_n$
with a basis $\{e_i\}$.  First we
construct an asymptotic version $Y_1 \in {\cal A}(X)$ such that
$E \in \{Y_1\}_n$. Let $\cal K$ be a family of  all $m$-tuples of
na\-tu\-ral numbers, for all  $m \le n$,  which are of the form
$K = \{p_1, p_1 p_2, \ldots, \prod_{i=1}^m p_i\}$,
where $p_1 < p_2 <\ldots < p_m$ are prime numbers.

For an arbitrary $N \in \NN$ sufficiently large,
let ${\cal K}_N \subset \cal K$  consists of
all $m$-tuples $K \in {\cal K}$, for all $m \le n$,  for which
$\prod_{i=1}^m p_i \le N$.
Observe that family  ${\cal K}_N $
has the property from Lemma~\ref{lemma_1}.
Therefore there exists an asymptotic  space $F_N \in \{X\}_N $
with a basis $\{f_i\}$ such that for any $m$-tuple
$K \in {\cal K}_N$ we have
$$
\{f_i\}_{i \in {K}} \stackrel{1}\sim \{e_l\}_{l=1}^{m}.
$$

Similarly as in~\ref{as_version},  we can then construct
an increasing sequence  of such spaces
$\ldots \subset F_N \subset F_{N+1}\subset \ldots$,
with bases
$\ldots \subset \{f_i\}_{i=1}^N
   \subset \{f_i\}_{i=1}^{N+1}\subset \ldots$,
each of them having the above structure
(because the restriction of $F_{N+1}$ to the first $N$
basis vectors has the same property). This sequence
defines an asymptotic version $Y_1$ and it can be checked
that $E \in \{Y_1\}_n$.

Given a finite number of asymptotic spaces
$\{E_{l}\}$,
we use Lemma~\ref{lemma_2} a finite number of times
to build an asymptotic version $\widetilde{Y}$
such that every $E_{l} \in \{\widetilde{Y}\}_n$.

The end of the argument is obvious: let $\ep_n \downarrow 0$
as $n \to \infty$. For every $n \in \NN$, let ${\cal T}_n$ be a finite
$\ep_n$-net in the set $\bigcup_{k \le n} \{X\}_k$
of all asymptotic spaces of dimension less than or equal to $n$.
Let $Y_n$ be an asymptotic version of $X$ which contains
all spaces from ${\cal T}_n$ as asymptotic spaces.
Use Lemma~\ref{lemma_3}  for the spaces $Y_n$ and the sets
$I_n = \{(2n+1) 2^t\}_{t=1}^\infty$. Resulting
asymptotic version $Y$ has the required property:
for every $n \in \NN$ and $E \in \{X\}_n$, there
is a sequence $\ep_k \to 0$ such that $Y$ has
permissible subspaces $F_k$ with the distance
$\di_b (F_k, E) \le 1 + \ep_k$. Thus $E \in \{Y\}_n$.
Therefore $ \{X\}_n \subset \{Y\}_n$; the converse inclusion
has been commented on before the statement of the theorem.
\qed


\section{Uniqueness of the asymptotic-$l_p$ structure}
\label{33}
\subsection{}
\label{thm_unique}
The following theorem has been already promised in~\ref{asymptotic_l_p_3}.

\begin{thm}
Let $X$ be a Banach space and consider the asymptotic structure 
on $X$ determined by the family $\cB^0 (X)$ of all 
finite-codimensional subspaces of $X$. 
Let $1 \le p < \infty$. 
Assume that there exists $C$ such that for every $n \in \NN$
and every $E \in \{X\}_n$, the Banach--Mazur
distance $\di(E, l_p^n) \le C$. Then $X$ is an 
asymptotic-$l_p$ space.
\end{thm}

\subsection{}
\label{proof_unique}
\proofsec
Let $\widetilde{Y}$  with the basis  $\{\widetilde{y}_i\}$
be a universal asymptotic version of $X$.
If $E \subset \widetilde{Y}$ is a finite-dimensional subspace
then for every $\ep >0$, $E$ is $(1+\ep)$-isomorphic to a subspace of
$Y_N = \spn [\widetilde{y}_i]_{i=1}^N$, which in turn
is $C$-isomorphic to $l_p^N$. Therefore 
$\widetilde{Y}$ is an ${\cal L}_p$-space. 
It is then well-known ([LR]) that   $\widetilde{Y}$  is isomorphic
to a  subspace  $Y$ of $L_p[0,1]$.
Let $\{y_i\}$ denotes the image of $\{\widetilde{y}_i\}$
by this isomorphism.
The asymptotic structure from $\widetilde{Y}$ induces 
an asymptotic structure on  $Y$, which we will now investigate.

Fix an arbitrary $E \in \{Y\}_n$ with the basis  $\{e_i\}$.
We will show that  there exists a permissible
$n$-tuple $\{z_i\} \sim \{e_i\}$ which is also equivalent
(up to some constant $D''$) to the unit vector basis in $l_p^n$.
This will  mean that $Y$ is an asymptotic-$l_p$ space, 
hence  by~\ref{isomorphic},   so is  $\widetilde{Y}$.
Since $\{X\}_n = \{\widetilde{Y}\}_n$ 
for every $n$, then $X$ itself will be  an 
asymptotic-$l_p$ space as well.

\subsubsection{}
\label{haar}
Consider the  Haar basis  in $L_p(0, 1)$.
It is easy to see from the vector game definition
that  $E$ can be realized
as successive blocks of this  basis.
Since the Haar basis  is unconditional for $1 < p < \infty$
(\cf~\eg~[LT.2]), then  in the case $1 < p < \infty$
the basis  $\{e_i\}$ is unconditional and 
$\mbox{unc}\{e_i\} \le K_p$,  where  $K_p$ depends on $p$ only.

This already completes the proof for $p=2$, since  
it is well-known that 
every unconditional basis in $l_2^n$ is equivalent 
to the standard unit vector  basis.
For $p \ne 2$ 
we proceed separately in  cases  $1 \le p < 2$ and  $ p > 2$.

\subsection{}
\label{p < 2}
Let  $1 \le p < 2$.
Fix $\ep >0$ and consider a vector game in $Y$ for $E$ and $\ep$.  Let
the first move of player \pS\   be $Y$ itself and let player 
\pV\   choose $z_1 \in S(Y)$. 
Considering appropriate choices for the second move of the
subspace player \pS, we obtain a  sequence 
$\{v_m\} \subset Y$ 
of second choices for \pV\   (with the first choice being always
$z_1$); we can  also  ensure that the $v_m$'s are  
successive blocks  of the Haar  basis.  

It is now convenient to describe the argument separately for
the reflexive case $ 1 < p < 2$ and for $p=1$. 

\subsubsection{}
\label{p < 2_beginning}
Let $ 1 < p < 2$.
Passing to a subsequence, we may assume that $\{v_m\} $
generates a spreading model (see~\ref{spreading_model}). 
It is known that the natural
basis of a spreading model of any sequence in  $L_p$
is symmetric, rather than spreading invariant as in general.
This is true in every stable Banach space, and  is  a direct 
consequence of the definition of stability 
[KM] (the reader not familiar with the notion of ultrafilters
may also consult [KM] p. 276). On the other hand,
$L_p$ is a stable space for $1 \le p < \infty$ ([KM]).
So from the definition of a spreading model, 
this means that for a given $n$, any  $n$-tuple
$v_{m_1},\ldots, v_{m_n}$, with $m_1< \ldots < m_n$ and 
$m_1$  large enough, 
forms   a finite almost symmetric  basis in its span.

On the other hand, given any weakly null sequence $\{v_m\}$ and  $M_1$,
by considering a suitable  subspace game we can choose
an $n$-tuple $v_{m_1},\ldots, v_{m_n}$ which is
2-equivalent to some  asymptotic  space and $m_1 \ge M_1$.
By the main assumption, the span of $\{v_{m_i}\}_{i=1}^n $
is $C$-isomorphic to $l_p^n$.

\subsubsection{}
\label{p = 1}
For  $p=1$  we need to be slightly more careful, because
the sequence $\{v_m\} $ is not weak null.
Still,  a finite analogue of the previous argument works here.
(It actually  does not require any assumptions on $p$ at all.)
 
Fix $N$ to be determined later, and note that 
in the definition of  $\{v_m\} $, by
using additionally a subspace game for $\{Y\}_N$ in $Y$, 
we can  also ensure that the vectors
$\{v_1, \ldots, v_N\}$  form a permissible $N$-tuple..
All infinite arguments from~\ref{p < 2_beginning} have finite analogues,
this follows from a standard compactness argument,
using  the stability of $L_1$ under ultraproducts
([DK]).
This means that given $n$,  there is $N$ such that
from every almost monotone normalized basic sequence 
$v_1, \ldots, v_N$, one can extract an almost 
symmetric  subsequence 
$v_{m_1},\ldots, v_{m_n}$ 
of length $n$. 
Since $\{v_1, \ldots, v_N\}$ is permissible, 
by~\ref{as_blocking}, 
this  subsequence can be assumed to be permissible
as well. In particular, as in~\ref{p < 2_beginning}, its span 
is $C$-isomorphic to $l_1^n$.

\subsubsection{}
\label{p < 2_jmst}
We go back to our more general assumption $1 \le p <2$.
Now we use [JMST] Theorem 1.5, which says that for 
$1 \le p \le \infty$, every $K$-symmetric basis in $l_p^n$ is 
$D'$-equivalent to the  standard unit vector basis 
$\{\mbox{e}_i\}$ in $l_p^n$, 
where $D' = D'(K)$ depends on $K$ only.
It follows that  there is $D = D(C)$ such that
$\{v_{m_i}\}_{i=1}^n $ is  $D$-equivalent 
to the basis $\{\mbox{e}_i\}_{i=1}^n$ in $l_p^n$.

\subsubsection{}
\label{p < 2_dor}
By Dor's result [D], which is valid for $1 \le p < \infty$,
for some $\delta = \delta (D) >0$,
there exist disjoint subsets $A_1,\ldots,A_n$ of $[0,1]$ such that
$$
 \int_{A_j} |v_{m_j}|^p \ge \delta \qquad \mbox{for  } j=1, \ldots, n.
$$
Since $\int |z_{1}|^p =1$,  
taking  $n$  sufficiently  large,  we get that at least 
one of the integrals 
$\int_{A_j} |z_{1}|^p $ is smaller than
$\delta / 4$; denote the corresponding set by $ A^{(2)}$ and   the
corresponding vector $v_{m_j}$ by $z_2$.

Passing  to a   sequence $\{w_m\}$ of possible third choices
for \pV\   in the vector game, with the first two choices  being  
$z_1, z_2$, and  repeating the argument we get
a set $A^{(3)}$ and a vector $z_3$ such that
$$
\int_{A^{(3)}} | z_3 |^p \ge \delta, \qquad 
\int_{A^{(3)}} | z_1 |^p < \delta /8, \qquad 
\int_{A^{(3)}} | z_2 |^p < \delta /4.
$$
By an obvious induction we get a permissible $n$-tuple 
$z_1, \ldots,  z_n$, 
$(1+\ep)$-equiv\-a\-lent to $\{e_i\}$   and 
disjoint sets $B_1, \ldots,  B_n$
such that 
$$
\int_{B_1} |z_1|^p \ge 1 - \delta/2,
\qquad 
\int_{B_i} |z_i|^p \ge \delta/2, 
\qquad 
\hbox{for } i>1.
$$
(We put 
$ B_1 =  ( \bigcup_{i=2}^n  A^{(i)})^c$
and $B_2 = A^{(2)} \backslash ( \bigcup_{i=3}^n  A^{(i)})$,
etc., ).
Recall that since $\{z_i\}$ are permissible, they are  unconditional.
Then the above condition implies a lower $l_p$-estimate:
for all $\{a_i\}$ we have
$$
\|\sum a_i z_i\|^p \sim \int (\sum |a_i z_i|^2)^{p/2}
\ge \sum_j \int_{B_j} |a_j z_j|^p \ge \delta /2 \sum |a_j|^p.
$$

\subsubsection{}
\label{p < 2_end}
The upper $l_p$-estimates are easy.
For   $1 < p < 2$ the estimate follows from the type $p$
and from the unconditionality of the basis 
$\{e_i\} $, obtained in~\ref{haar}. 
For $p=1$,  we use the triangle inequality.

Thus $\{z_i\}_{i=1}^n$ is $D_p''$-equivalent to the unit
vector basis in $l_p^n$, as required, where the constant 
$D_p''$ depends on $C$ and on $1 \le p < 2 $.

\subsection{}
\label{p > 2}
Let $ p > 2$.
We use Kadec--Pe{\l}czy\'{n}ski approach 
(\cf~[LT.2], 1.c.8).
For $x \in L_p$  and $\delta >0$, set
$\sigma (x, \delta) = \{t\in [0, 1] \, 
\mid\, |x(t)| \ge \delta \|x\|\}$,
and let 
$M(\delta) = \{x \,\mid\,  \mu (\sigma (x, \delta))< \delta\}$.

We start with a couple of general remarks which  can be proved
by standard well-known arguments.

\subsubsection{}
\label{p > 2_gen}
Recall that if  a sequence of functions 
$\{w_m\}$  is $K$-unconditional and it belongs to  $M(\delta)$,
for some $\delta >0$, then  
$\{w_m\}$  satisfies a lower $l_2$ estimate
(with a constant depending on $K$ and $\delta$)
(\cf~\eg~[LT.2], 1.c.10). If $p >2$, combining this
with the type 2 of the space $L_p$  we get that 
$\{w_m\}$  is equivalent to the unit vector basis in $l_2$.

\subsubsection{}
\label{p > 2_induc}
Consider a sequence $\{w_m\} $
such that $w_m \not\in M(2^{-m-2})$
and let  $\eta_m = \sigma (w_m, 2^{-m-2})$, for 
$m=1, 2, \ldots$.  
Given $z_1, \ldots, z_k$ in $L_p$, there exists
$m_0$ such that
$\int_{\eta_{m_0}}   |z_i|^p < 2^{-k-2}$
for $i=1, \ldots, k$.

\subsubsection{}
\label{p > 2_game}
Recall that $E \in \{Y\}_n$ was an arbitrary asymptotic space with a
basis $\{e_i\}$ and consider the same games for $E$ as in~\ref{p < 2}.
Let us outline an inductive argument. Let $0 \le k < n$ and
assume that $z_1, \ldots, z_k$
have been already defined as possible choices for the first $k$ moves
of player \pV\   (in a vector game for $E$). With these vectors fixed,
consider a $w$-null sequence $\{w_m\}$ of possible choices in the
$(k+1)$th move for \pV.  Using~\ref{p > 2_gen} and our main
isomorphism  assumption,
we conclude  that there is no $\delta$ such that 
$\{w_m\} \subset M(\delta)$. Passing to a subsequence 
we may therefore assume that 
$w_m \not\in M(2^{-m-2})$, for $m=1, 2, \ldots$.
Let $m_0$ be as in~\ref{p > 2_induc}, denote $w_{m_0}$ by
$z_{k+1} $ and set $\sigma_{k+1}= \eta_{m_0}$.

Proceeding this way we get a permissible $n$-tuple $\{z_i\}$,
$(1+\ep)$-isomorphic to $\{e_i\}$, and subsets $\sigma_i$
of $[0,1]$, such that  for every $i=1, \ldots, n$ we have
$$
\int_{\sigma_k}   |z_i|^p < 2^{-k-2}
\qquad \mbox{for}\quad i < k \le n.
$$

Obviously, for  $i=1, \ldots, n$ we have
$\int_{\sigma_i^c}   |z_i|^p 
<  2^{(-i-2)p}$. Thus, setting
$B_i = \sigma_i \backslash \bigcup_{k>i} \sigma_k$
we get
$$
\int_{ B_i^c}   |z_i|^p <  2^{(-i-2)p} + \sum_{k >i} 2^{-k-2}
   <  2^{-i-1}
\qquad \mbox{for}\quad i=1, \ldots, n.
$$

Since  $\|z_i\| =1$ for  $i=1, \ldots, n$, then $\{z_i\}$
are equivalent (up to a universal constant) to the unit vector
basis in $l_p^n$.  As already indicated  at the end
of~\ref{proof_unique}, this completes the proof of the theorem.
\qed


\section{Duality for asymptotic-$l_p$ spaces} 
\label{44}
\subsection{}
\label{min_system}
A minimal system  in a Banach space $X$ is a sequence $\{u_i\}$
such that there exists a sequence $\{u_i^*\}$ in $X^*$
so that $\{u_i, u_i^*\}$ is a biorthogonal system.
Systems considered here will be always  fundamental and total,
in particular, $X = \overline{\spn}\{u_i\}$.
Some more information, and in particular
classical definitions
of shrinking and boundedly complete minimal systems,
can be found  \eg~in [LT.1], I.f. 
Let us just recall that a space $X$ is
reflexive if and only if every minimal system
in $X$ is both shrinking and boundedly complete.
The reader who is not familiar with minimal systems
may just think about a basis in $X$.

Recall that $\cB^t(X)$  denotes
the family of all tail subspaces.

\subsubsection{}
\label{mi_sh_2.1}
Let us recall the following known fact
([Mi], also [MiS], Proposition 2.1). In
presence of a basis in $X$ this fact is obvious and 
does not require the shrinking   assumption.

\begin{lemma}
Let $(Y, \|\cdot\|_Y)$ be a Banach space with a shrinking 
minimal system.  There exists an equivalent norm  $\|\cdot\|$
on $Y$ such that $\|x\|_Y \le \|x\| \le 2 \|x\|_Y$ for all
$x \in Y$ and that for every $\delta >0$ and every tail subspace
$\widetilde{Z} \in \cB^t(Y^*)$ there exists a tail subspace
$\widetilde{Y} \in \cB^t(Y)$  such that for every 
$x \in S(\widetilde{Y})$ there is $f \in S(\widetilde{Z})$
with $f(x) \ge 1 - \delta$.
\end{lemma}

\subsection{}
\label{shrink_bdd_comp}
Let $X$ be an asymptotic-$l_p$
space (with respect to  the family $\cB^0(X)$)
and let  $\{u_i\}$ be a  minimal system  in  $X$.
If  $ 1 < p \le \infty$, the system  is shrinking. 
If  $ 1 \le p  < \infty$, the system
is boundedly complete.

Assume to the contrary that  $\{u_i\}$ is not shrinking,
\ie\  $X^* \ne \overline{\spn}\{u_i^*\}$.
Fix $n\in \NN$ to be defined later.
There exists $x^* \in X^*$ with $\|x^*\|= 1$,
for which one can construct a permissible $n$-tuple  
$\{x_i\}$  in $X$ (for an arbitrary $\ep >0$),
such that  $|x^*(x_i)| >\delta$ for $i=1, \ldots, n$,
where  $\delta >0$ is a universal constant.
Then 
$$
C (1 + \ep) n^{1/p} \ge
\Bigl\|\sum_{i=1}^{n} x_i\Bigr\|
\ge \Bigl|x^*\Bigl(\sum_{i=1}^{n} x_i\Bigr)\Bigr|\ge n \delta,
$$
and, if $p >1$, this is a contradiction  for $n$ large enough.

Assume that  $\{u_i\}$ is not boundedly complete.
For every $n \in \NN$, there exists a permissible
(normalized) $n$-tuple $\{x_i\}$ such that
$\sup_n \| \sum _{i=1}^n x_i\|= M <\infty$,
where $M$ is a universal constant.
On the other hand
$ \| \sum _{i=1}^{n} x_i\| \ge (1/C)\, n^{1/p}$,
which is a contradiction,
if $p < \infty$.

In particular, for $ 1 < p < \infty$, an asymptotic-$l_p$
space is reflexive. Note however that $l_1$-  and $l_\infty$-spaces
may be reflexive as well; such examples are given by
the Tsirelson space  $T_{(1)}$  and its dual $T_{(1)}^*$ 
(\cf~\eg~[CS]).

\subsection{}
\label{duality}
\begin{thm}
Let $ 1 \le p \le \infty$ and let $X$ be an asymptotic-$l_p$
space which is reflexive. 
Then $X^*$ is an asymptotic-$l_{p'}$, where $ 1/p + 1/p' =1$
(with the standard convention for $p=1$ and $p= \infty$).
\end{thm}

\subsubsection{}
\label{extra_comments}
Let $\{u_i\}$ be a minimal system in  $X$.
By~\ref{isomorphic} we may assume,
without loss of generality,  that
the norm in $X$ satisfies the conclusion of
Lemma~\ref{mi_sh_2.1}.
Moreover,  the asymptotic structures of $X$ and of $X^*$
are   determined by the families $\cB^t(X)$
and $\cB^t(X^*)$   associated to $\{u_i\}$  and to $\{u_i^*\}$,
respectively.

To make the statements below more intuitively clear
and to avoid tiresome repetitions,
let us recall (\cf~\ref{vector_game_1}
and~\ref{subspace_game_warning}) that
if $\ep >0$ is fixed,
then  an $n$-tuple in $X$ (resp.\ in $X^*$)
is permissible,  if it is $(1+\ep)$-equivalent
to the natural basis in an asymptotic
space from $\{X\}_n$ (resp. $\{X^*\}_n$).
In particular an $n$-tuple is permissible
if it is obtained as a result of a subspace game
in $X$ (resp.\ in $X^*$),
assuming that  player \pS\ followed his winning
strategy for $\{X\}_n$ (resp. $\{X^*\}_n$) and $\ep$.

\subsection{}
\label{easy_est}
An asymptotic lower $l_{p'}$ estimate in $X^*$
is based on the following  lemma.

\begin{lemma}
Let $Y$ be a Banach space with a shrinking minimal system.
Let  $\{e_i\} \in \{Y\}_n$ be an asymptotic $n$-tuple
and let $\ep >0$.
There exist a permissible $n$-tuple $\{z_i\}$ in $Y$
satisfying $\{z_i\} \stackrel{1+\ep} \sim   \{e_i\}$,
and a permissible  $n$-tuple $\{g_i\} \subset S(Y^*)$ in $Y^*$,
such that $g_i (z_i) \ge 1 - \ep$ for $i=1, \ldots, n$
and  $g_i (z_j) =0$ if $i \ne j$.
\end{lemma}

The proof of the lemma
requires  an asymptotic game in $Y$, which
combines strategies for two simultaneous
games: a winning strategy for \pV\ in a vector game in $Y$ and
a winning strategy for \pS\ in a subspace game in $Y^*$.
The latter strategy  ensures permissibility in $Y^*$ and
it   determines choices of subspaces in $Y$ via~\ref{mi_sh_2.1}
(\cf~the proof of  Lemma~\ref{global_est} below).
We leave it for the reader.

Now the proof of the lower $l_{p'}$ estimate
follows   a standard argument.
Given an asymptotic $n$-tuple $\{e_i\}$ in $X^*$
and $\ep >0$, let   $\{z_i\}$  in $X^*$ and $\{g_i\}$
in $X$ be as in the lemma. For any scalar $n$-tuple
$a = \{a_i\}$, pick $b= \{b_i\}$ with $\|b\|_p= 1$ such that
$\sum_i a_i b_i = \|a\|_{p'}$.
Then
$$
(1 - \ep) (\sum_i |a_i|^{p'})^{1/{p'}} \le
(\sum_i b_i g_i) (\sum_i a_i z_i)\le
C (1+\ep) \|\sum_i a_i e_i\|,
$$
as required.

\subsection{}
\label{global_est}
An asymptotic upper $l_{p'}$ estimate in $X^*$
is based on the following  reformulation
in our context of Theorem 2.2 from [MiS].

\begin{lemma}
Let $Y$ be a Banach space with a shrinking minimal system.
Let  $\{e_i\} \in \{Y\}_n$ be an asymptotic $n$-tuple,
let $\{a_i\}$ be an arbitrary scalar  sequence
and  let $\ep >0$.
There exist a permissible $n$-tuple $\{y_i\}$ in $Y$
satisfying
$\{y_i\} \stackrel{1+\ep} \sim   \{e_i\}$,
and a permissible  $n$-tuple $\{g_i\}$ in $Y^*$,
and a sequence of  scalars $\{b_i\}$,
such that $g_i (y_j)= 0$ if $i \ne j$ and
$$
\Bigl(\sum_{i=1}^n b_i g_i\Bigr)\,\Bigl(\sum_{j=1}^n a_j y_j\Bigr)
\ge (1-\ep) \, \Bigl\|\sum_{i=1}^n b_i g_i\Bigr\|\
\Bigl\|\sum_{j=1}^n a_j y_j\Bigr\|.
$$
\end{lemma}

\subsubsection{}
\label{global_est_proof}
This result  is based on an argument which might be useful
in other context;  for the reader convenience we outline the proof.

\smallskip
\proof
We provide a complete argument for $n=2$, with few comments
concerning the general  case.
Let $Z = Y^*$ and fix $\delta >0$ to be defined later.
Consider a subspace game  in $Z$ for $\{Z\}_2$ and $\ep$.
We name the players of this game by \pS$^*$ and \pV$^*$
respectively.  Let $Z_1   \in \cB^t(Z)$
be a tail subspace chosen by   \pS$^*$
in the first move. Let ${Y_1} \in \cB^t(Y)$ be a corresponding
subspace  (for $\delta$), as in  Lemma~\ref{mi_sh_2.1}.
Now consider a vector game in $Y$ for
$\{e_i\}$ and $\ep$, with ${Y_1}$  being
the first choice of  player \pS.
Let  player \pV\    choose $y \in S(Y_1)$.
Considering  appropriate choices for the
second  move of  the subspace player \pS,
we obtain a  sequence   of successive blocks
$y_1 < y_2 < \ldots$ of second choices for \pV\
(with the first choice always being  $y$).
(If  $n > 2$, then with  a fixed $m \in \NN$
let  $\{y_{m,l}\} \subset S(Y_1) $ be a
sequence of  successive blocks,
each of which  could be  picked by   \pV\   in his third move,
in the game in which his first two moves were
$y$ and $y_m$.  And so on.)

Fix  $m$.  Then  $\{y, y_m \}\stackrel{1+\ep} \sim \{e_1, e_2\}$.
Let $w_{ m}= a_1 y + a_2 y_m $.
Let $f_{m} \in S(Z_1)$ be a functional  norming
$w_{ m}$ up to $\delta$,  as  in
Lemma~\ref{mi_sh_2.1}.
We will  show that there is $\mu \in \NN$ such that $f_\mu$
can be approximated (up to $3 \delta$) by a functional
of a form $h= b_1 g_1 + b_2 g_2$,  with  $\{g_i\}$
permissible and satisfying the required
biorthogonality condition.
In particular, $h$ will norm $w_\mu$ up
to $4\delta$,  which will give the conclusion
by setting $\delta = \ep / 4$.

Let $f$ be a $w^*$-cluster point of
$\{f_{m}\}_m$.  Then $f \in Z_1$.
Let $h_1 \in Z_1$ be  finitely supported
such that $\|h_1 - f\| < \delta$.
Let $g_1 = h_1/ \|h_1\|$ and consider
this $g_1$ as a choice for the vector player \pV$^*$
in the subspace game in $Z$,
so that \pS$^*$ chooses $Z_1$ and \pV$^*$ chooses $g_1$.
Let $Z_2 \in \cB^t(Z)$
be a subspace   picked by  \pS$^*$
in his second move. Then $Z_2$ is  the $k$th tail subspace,
for some $k \in \NN$
and we may assume without loss of generality that
$k > \max (\supp (g_1)\cup \supp (y))$.
Let $Q_{k}$ denote the canonical projection
in $Z$ onto $\spn \{u_i^*\}_{i \le {k}}$,
so that in particular $Q_{k} h_1 = h_1$.
Pick  $\mu \in \NN$  such that
$\|Q_{k} f_{\mu} - h_1 \| \le \|Q_{k} (f_{\mu} - f) \|
   +  \|Q_{k} (f - h_1) \| < 2 \delta$
and that $\min \supp (y_\mu) >k $.
Then $(I - Q_{k}) f_{\mu} \in Z_2$
and pick finitely supported  $h_2 \in Z_2$  such that
$\|(I - Q_{k}) f_{\mu}- h_2\|< \delta$. Set $g_2 = h_2 / \|h_2\|$.
Then $\{y, y_\mu\}$ is the required permissible couple in $Y$.
Note that $g_1$ and $y_\mu$ are disjointly supported,
and so are $y$ and $g_2$. Thus $g_1(y_\mu) = g_2 (y) =0$.
Also,  the functional $h = h_1 + h_2$   approximates
$f_\mu$ up to $3 \delta$,  as promised.
Finally, since $Z_2$ was  a second choice of \pS$^*$ and
$g_2 \in Z_2$, then $\{g_1, g_2\}$ is a permissible
couple in $X^*$ and of course, $h = b_1 g_1 + b_2 g_2$,
for suitable scalars $b_1, b_2$.
(If $n >2$,  consider  $g_2$ as a second choice for
\pV$^*$,  and let  $Z_3$ be a subspace picked
by \pS$^*$, which starts after $g_2$ and $y_\mu$
and then repeat the argument.)
\qed

\subsubsection{}
\label{upper_l_p'}
Again, the proof of an upper $l_{p'}$-estimate is completely
standard. Given an asymptotic $n$-tuple $\{e_i\}$ in $X^*$,
scalars  $ \{a_i\}$ and $\ep >0$,
apply the lemma for $Y = X^*$ to get
$\{y_i\}$  in $X^*$ and $\{g_i\}$ in $X$ and scalars $\{b_i\}$,
with the additional  normalization $\|\sum_i b_i g_i\|=1$.
Then $(\sum |b_i|^p)^{1/p}\le C$. Thus
$$
(1-\ep) \, \Bigl\|\sum_{j=1}^n a_j y_j\Bigr\|
\le
\Bigl(\sum_{i=1}^n b_i g_i\Bigr)\,\Bigl(\sum_{j=1}^n a_j y_j\Bigr)
=  \sum_{i=1}^n a_i  b_i \le
  C\, \Bigl(\sum_{i=1}^n |a_i|^{p'}\Bigr)^{1/p'},
$$
as required.
Combined with~\ref{easy_est}, this concludes the proof of
Theorem~\ref{duality}.
\qed


\section{Complemented permissible subspaces} 
\label{55}

It is well-known and easy to see that every block subspace
of $l_p$ is complemented; the same is true for Tsirelson spaces
$T_{(p)}$, although in this case it is much more difficult
to prove (here $1 \le p <  \infty$)(\cf~[CS]). 
To get a related complementation property which
actually characterizes spaces $l_p$ or $c_0$, one needs to add 
an unconditionality assumption and to consider all permutations
of a given basis ([LT.3], \cf~{\em also} [LT.1] 2.a.10).
In the asymptotic setting  the situation is 
more natural and elegant, and a natural complementation
condition fully characterizes asym\-pto\-tic-$l_p$ spaces.

\subsection{}
\label{asymp_properties}
We start  by describing  few more 
asymptotic notions. 
Let  ${\cal P}$ be  a property of finite-dimensional 
subspaces of  a given Banach space.
\begin{dfn}
We say that ${\cal P}$ is satisfied by 
{\it permissible subspaces of $X$ far enough}, 
if for every $n\in \NN$ and $\ep > 0$,
the subspace player
\pS\   has a winning strategy in a subspace game for 
$\{ X\}_n$ and $\ep$ such that arbitrary $n$-tuple $\{x_i\}$
resulting  from the game  spans a  subspace with 
property ${\cal P}$. (This subspace is automatically
permissible, since the strategy is winning for $\{X\}_n$.)
\end{dfn}

We have a similar definition if  ${\cal P}_n$ is a 
property of $n$-dimensional  subspaces of $X$, 
with $n\in \NN$  fixed. 

Any strategy for \pS\ as above will be called  
a  ${\cal P}$-strategy.

\subsubsection{}
\label{permis_property_tree}
Recall our  intuition of a tree-like structure of subspaces and
vectors, as in~\ref{subspace_game_tree}.
Then ${\cal P}$ is satisfied by {\it permissible subspaces of $X$ 
far enough}, if   and only if for an
arbitrary $n \in \NN$ and $\ep >0$, by pushing subspaces {\it far
enough\/} along $\cB(X)$ player \pS\ can ensure that the 
subspaces spaned  by all resulting $n$-tuples
are not only permissible but  they  also have property
${\cal P}$.

\subsubsection{}
\label{permis_property_V}
Assume that ${\cal P}$ is
satisfied by  permissible subspaces of $X$ far enough
and let  $n \in \NN$ and $\ep >0$.
By combining the strategy for the subspace player \pS\  
discussed in~\ref{as_blocking} with a $\cal P$-strategy,
and using  filtration condition~\ref{filtration},
we obtain a strategy for \pS\  such that
arbitrary normalized successive blocks 
of any  $n$-tuple $\{x_i\}$
resulting  from the game, are permissible and their span has
property ${\cal P}$.

Of course  the set of all $n$-tuples resulting from the 
game above represents all spaces  from $\{X\}_n$.
In other words, for any $E \in \{X\}_n$ there is 
$\{x_i\}$ as above, $(1+\ep)$-equivalent to the basis in $E$.
Indeed, we could   appropriately instruct  player \pV\ 
to achieve this $E$, up to $ 1+ \ep$.

\subsubsection{}
\label{permis_projection}
Let $X$ be a Banach space with a minimal system $\{u_i\}$.
Let $Y = \spn [y_i]$ be a 
subspace of $X$. A projection $P: X \to Y$ is 
called    $\{u_i\}$-{\it permissible\/}  
(or just  {\it permissible}, if the system 
$\{u_i\}$ is understood from the context)
if $P$ can be written as
$  P = \sum_i g_i \otimes y_i $,
with $g_i \in X^*$ finitely supported
and $\max \supp(g_i) < \min \supp (g_{i+1})$,
for $i=1, 2, \ldots$.

\subsection{}
\label{asympt_l_p_XX}
The duality theorem~\ref{duality} implies 
(and in fact is equivalent to)
a complementation property of   asymptotic-$l_p$ spaces.
\begin{cor}
Let  $X$ be an asymptotic-$l_p$ space
for some $1 < p < \infty$.
Then there is $D$ such that  permissible subspaces 
of $X$ far enough are $D$-com\-ple\-mented by means of 
permissible projections.
\end{cor}
\proof
Let $\{u_i\}$ be a minimal system in  $X$, and without loss of generality
let us make all the assumptions as in~\ref{extra_comments}.
Let $n \in \NN$ and  $\ep >0$. Player \pS\ has a strategy
in an asymptotic game in $X$ such that 
if $\{x_i\}$ is a resulting permissible $n$-tuple  in $X$,
then there exists 
a permissible  $n$-tuple $\{g_i\} \subset S(X^*)$,
with $\max \supp(g_i) < \min \supp (g_{i+1})$, 
such that $g_i (x_i) \ge 1 - \ep $ for $i=1, \ldots, n$
and  $g_i (x_j) =0$ if $i \ne j$.
Indeed, the strategy for \pS\ is  essentially the same
as in Lemma~\ref{easy_est}, with an additional requirement for
successiveness of the $g_i$'s.
This property formally implies the existence of a required 
permissible projection onto  $\spn [x_i]$.

Let $  P = \sum_i g_i \otimes x_i $. Clearly, $P$ is a 
permissible projection onto $\spn [x_i]$.
Fix an arbitrary vector $x \in X$ and pick
scalars $\{b_i\}$  such that 
$\sum_i g_i(x) b_i = (\sum_i |g_i(x)|^p)^{1/p}$
and $\sum_i |b_i|^{p'} =1$. Since $X$ is
asymptotic-$l_p$ space and, by
Theorem~\ref{duality},  $X^*$ is asymptotic-$l_{p'}$ space
(with a constant $C$), then 
\begin{eqnarray*}
  \|Px\| &= & \|\sum_i g_i(x) x_i\|\le C  (\sum_i |g_i(x)|^p)^{1/p}\\
   &=&  C \sum_i g_i(x) b_i \le C \|\sum_i b_i g_i\| \le C^2.
\end{eqnarray*}
Thus $\|P\| \le C^2$.
\qed

\subsection{}
\label{complemented}
For spaces with basis  the converse is true.

\begin{thm}
Let $X$ be a Banach space with a basis. Assume that 
there exists a constant $C$ such that 
permissible subspaces of $X$ far enough are $C$-complemented 
by means of  permissible projections. 
Then  $X$ is an asymptotic-$l_p$ space
for some $1 \le p \le \infty$.
\end{thm}

The asymptotic structure in $X$ may be naturally taken 
with respect to either family $\cB^0(X)$ or $\cB^t(X)$.
Then the conclusion of the theorem relates to the same structure.

\subsection{}
\label{comments_complem}
Before we pass to the proof of the theorem, let us make 
some comments.

\subsubsection{}
\label{comm_1}
The argument below shows that if the basis in $X$ is unconditional
then the assumption that projections are permissible
can be dropped.

\subsubsection{}
\label{comm_2}
For arbitrary Banach spaces we have

\begin{cor}
Let $X$ be a Banach space. Assume that a universal 
asymptotic space
$Y \in {\cal A}(X)$ has the property that   
there exists a constant $C$ such that 
permissible subspaces of $Y$ far enough are $C$-complemented
(in $Y$) by means of  permissible projections. 
Then  $X$ is an asymptotic-$l_p$ space
for some $1 \le p \le \infty$.
\end{cor}

This corollary follows immediately by applying 
Theorem~\ref{complemented}  to $Y$.

\subsection{}
\label{proof_complem}
The argument below  is an asymptotic analogue of the original proof 
as presented \eg~in [LT.1] 2.a.10.

\proof
Let $n \in \NN$ and  let $ \{v_i\} \in \{X\}_n$
be an asymptotic $n$-tuple.
Fix Krivine's  $p \in [1, \infty]$, as in~\ref{krivine}.
Let $I_1 = \{k(n+1)+1 \, \mid \, k=0, \ldots, n-1\}$
and let $I_2 = \{1, \ldots, (n+1)^2\} \backslash I_1$.
By~\ref{mixtures}, there exists an  asymptotic 
$(n+1)^2$-tuple $\{f_l\}$ such that
$\{f_l\}_{l \in {I_1}} \stackrel{1} \sim \{v_i\}$
and  $\{f_l\}_{l \in {I_2}} \stackrel{1} \sim \{\mbox{e}_k\}$,
where  $\{\mbox{e}_k\}$ is the unit vector basis in 
$l_p^{n(n+1)}$.

Now fix $\ep >0$ and let $\{u_l\}$ be a permissible 
$(n+1)^2$-tuple of successive blocks of the basis in $X$,
$(1+\ep)$-equivalent to $\{f_l\}$ and 
such that    all subspaces spaned by  successive blocks
of $\{u_l\}$  admit permissible  projections of norm $\le C$.
This is possible by  the final comment in~\ref{permis_property_V}.
Set $F = \spn [u_l]_{l=1}^{n(n+1)}$ and  
$E = \spn [u_l]_{l \in {I_1}}$ and relabel
the basis in $E$ by $\{x_i\}_{i=1}^n$.

\subsubsection{}
\label{2}
By the assumption, there exists a projection 
$Q: F \to E$ with $\|Q\| \le C$. 
Since $\codim \ker Q = n$, for 
every $j= 1, \ldots, n$, we can find vectors
$e_j \in \ker Q \cap \spn [ u_{(j-1)(n+1) +1}, \ldots,  u_{j(n+1)}]  $
with $\|e_j \| =1$.
Thus we have successive blocks of the basis
$ x_1, e_1, x_2, \ldots, x_n,  e_n $
and we denote their span by  $Z$.
Of course, $\{e_j\} \stackrel{1+\ep}\sim \{\mbox{e}_i\}$, and
it  will cause no confusion to write $ l_p^{n}$ for $\spn [e_j]$.
By the construction,
$Z$ is  a permissible $2n$-dimensional subspace of $X$
and $Z = E \oplus l_p^n$;
the natural projection $Q$ on the first
coordinate has norm $\le C$
(hence the norm of  the projection on
the second coordinate is $\le C+1$).

\subsubsection{}
\label{3}
Fix $\la >0$ and let
$$
G = \spn[x_1 + \la e_1, x_2 + \la e_2, \ldots, x_n + \la e_n]
 \subset Z.
$$
Since $G$ is a block  subspace of $Z$, there is  a permissible
projection $P: Z \to G$ onto $G$ with $\|P\| \le C$.
The form of $G$ implies that   $P$ written  in the $C$-direct
sum  decomposition of $Z$ has a  matrix
of the form
$$
P_{|Z} = \left[
      \begin{array}{cc}
 A & (1/\la) B\\
\la A & B
\end{array}
\right]
$$
In other words, writting
$ P = \sum_i z_i^* \otimes (x_i + \la e_i)$, as
in~\ref{permis_projection}, we have
$$
a_{i,j} = z_i^* (x_j) \mbox{\ \  and \ \ } b_{i,j} = \la z_i^* (e_j)
\qquad \mbox{for } i, j = 1, \ldots, n.
$$

\subsubsection{}
\label{4}
Since $P$ is a projection, we have $A + B = I$; that is,
$$
a_{i,j}+  b_{i,j} = z_i^* (x_j + \la e_j) = \delta_{ij}
\qquad \mbox{for } i, j = 1, \ldots, n.
$$
Since $P$ is permissible, we have
$\max \supp z_i^* < \min \supp z_{i+1}^*$,
for all $i=1, \ldots, n-1$.
Finally, the form of $P$ implies the norm estimates:
$$
  \|A: E \to l_p^n \| \le (1/ \la)C(C+1)  \qquad\mbox{and}\qquad
  \|B: l_p^n \to E \| \le \la C^2.
$$

Since $\supp z_i^* \cap \supp (x_i + \la e_i) \ne \emptyset$,
and $ \max \supp z_{i-1}^* <  \min \supp z_{i}^*$
and $\max \supp z_{i}^*  < \min \supp z_{i+1}^* $,  then
$\supp  z_i^* \cap \supp x_{j} = \emptyset$,
if $|i - j|>1$ and  $i=1, \ldots, n$.
In particular, $a_{i, j}=0$ if $|i - j|>1$ and
$i=1, \ldots, n$.
Similarly,  $b_{i, j}=0$ if $|i - j|>1$ and
$i=1, \ldots, n$.
So for any $\la >0$, the matrices of  operators
$A$ and $B$  are tri-diagonal.

\subsubsection{}
\label{5}
Let $\la = 1/4 C^2$. Then, by~\ref{4},
$|b_{i,j}| \le   \|B: l_p^n \to E \| \le 1/4 $,
for $i,j =1, \ldots, n$.
Since the matrix of $B$ is tri-diagonal,
$\|B: l_p^n \to l_p^n \|\le 3 \max_{i,j} |b_{i,j}|\le 3/4$.
Since $I - A = B$, this implies that
$A$ is invertible on $l_p^n$ and
$$
\|A^{-1}: l_p^n \to l_p^n \|\le 4.
$$
Combining with norm estimates from~\ref{4} we get
$$
\|I: E \to l_p^n \| \le
  \|A: E \to l_p^n \|\  \|A^{-1}: l_p^n \to l_p^n \|\le 16 C^3 (C+1).
$$
This means that the vectors $\{x_i\}$ satisfy the lower
$l_p$-estimate with the constant $C'= 16 C^3(C+1)$.

\subsubsection{}
\label{6}
Let $\la = 4 C(C+1)$.  By~\ref{4}, $|a_{i,j}|\le 1/4$,
for $i,j =1, \ldots, n$, hence
$\|A: l_p^n \to l_p^n \|\le 3/4$.
Thus $\|B^{-1}: l_p^n \to l_p^n \|\le 4$, and hence
$$
\|I:  l_p^n  \to E \| \le
 \|B^{-1}: l_p^n \to l_p^n \|\ \|B:  l_p^n  \to E  \|
  \le 16 C^3(C+1) = C'.
$$
It follows that  the vectors $\{x_i\}$ satisfy the upper
$l_p$-estimate with the constant $C'$.
Thus, $\{x_i\} \stackrel{{C'}^2}\sim \{\mbox{e}_i\}$.
By the construction at the beginning of the proof, the same
holds for $\{v_i\} \in \{X\}_n$, hence
$X$ is an asymptotic-$l_p$.
\qed


\address

\begin{thebibliography}{JMST}
%
\bibitem[BL]{}  {\sc Beauzamy, B.} \& {\sc Laprest\'{e}, J.-T.},
  ``Mod{\`e}les {\'E}tal{\'e}s des Espaces de Banach'',
  Hermann, 1984.
\bibitem[BS]{}  {\sc Brunel, A.} \& {\sc Sucheston, L.},
  On B-convex Banach spaces, 
  {\sl Math.~Syst.~Th.}, 7 (1974), 294--299.
\bibitem[C]{}   {\sc Casazza, P.~G.},
  Some questions arising from the homogeneous Banach  space problem,
  in ``Banach Spaces'', {\sl Contemp. Math.}, 144 (1993), 35--52.
\bibitem[CS]{}  {\sc Casazza, P.~G.} \& {\sc Shura, T.},
  ``Tsirelson's Space'',
  Lecture Notes in Math., 1363, Springer Verlag, 1989. 
\bibitem[DK]{}  {\sc Dacunha-Castelle, D.} \& {\sc Krivine, J.-L.},
  Application des ultraproduits \`a l'\'etude des espaces 
  et des alg\`ebres de Banach, 
  {\sl Studia Math.}, 41 (1972), 315--334.
\bibitem[D]{}   {\sc Dor, L.~E.},
  On projections in $L_1$, {\sl Ann.\ of Math.}, 102 (1975), 463--474.
\bibitem[G.1]{} {\sc Gowers, W.~T.},
  A new dichotomy for Banach spaces, preprint
\bibitem[G.2]{} {\sc Gowers, W.~T.},
  A hereditary indecomposable space with an asymptotically unconditional
  basis, preprint.
\bibitem[G.3]{} {\sc Gowers, W.~T.},
  A Banach space not containing $l_1$ or $c_0$ or a reflexive subspace,
  preprint.
\bibitem[GM]{} {\sc Gowers, W.~T.} \& {\sc Maurey, B.},
  The unconditional basic sequence problem, 
  {\sl Journal of AMS}, 6 (1993),  851--874.
\bibitem[JMST]{}{\sc Johnson, W.~B.} \& {\sc Maurey, B.} \&
                {\sc Schechtman, G.} \& {\sc Tzafriri, L.},
  ``Symmetric Structures in Banach Spaces'',
  Memoirs of the AMS, no.~217, vol.~19, AMS, 1979.
\bibitem[K]{}   {\sc Krivine, J.-L.}, 
  Sous-espaces de dimension finie des espaces de Banach r{\'e}ticul{\'e}s,
  {\sl Ann.\ of Math.}, 104 (1976), 1-29.
\bibitem[KM]{}  {\sc Krivine, J.-L.} \& {\sc Maurey, B.},
  Espaces de Banach stables, {\sl Israel J. Math.}, 39 (1981), 273--295.
\bibitem[LR]{}{\sc Lindenstrauss, J.} \& {\sc Rosenthal, H.~P.},
  The ${\cal L}_p$ spaces, {\sl Israel J. Math.}, 7 (1969), 325--349.
\bibitem[LT.1]{}{\sc Lindenstrauss, J.} \& {\sc Tzafriri, L.},
  ``Classical Banach Spaces I, Sequence Spaces'',
  Springer Verlag, 1977.
\bibitem[LT.2]{}{\sc Lindenstrauss, J.} \& {\sc Tzafriri, L.},
  ``Classical Banach Spaces II'',
  Springer Verlag, 1979.
\bibitem[LT.3]{}{\sc Lindenstrauss, J.} \& {\sc Tzafriri, L.},
  On the complemented subspaces problem,
  {\sl Israel J. Math.}, 11 (1971), 263--269.
\bibitem[Ma]{}  {\sc Maurey, B.}, A remark on distortions, preprint.
\bibitem[Mi.1]{}  {\sc Milman, V.~D.},
  The geometric theory of Banach spaces, Part II, 
  {\sl Usp. Mat. Nauk}, 26 (1971), 73--149 (in Russian),
  (English translation: Russian Math.~Surveys 26 (1971),
  79--163).
\bibitem[Mi.2]{}  {\sc Milman, V.~D.},
  Spectrum of continuous bounded functions on the init shere of a Banach
  space, {\sl Funct.~Anal.~Appl.},  3 (1969), 67--79.
\bibitem[MiS]{}{\sc Milman, V.~D.} \& {\sc Sharir, M.},
  Shrinking minimal systems and complementation of $l_p^n$-spaces
  in reflexive Banach spaces, 
  {\sl Proc.~London Math.~Soc.}, 39 (1979), 1--29.
\bibitem[MiSch]{}{\sc Milman, V.~D.} \& {\sc Schechtman, G.},
  ``Asymptotic Theory of Finite Dimensional Normed Spaces'',
  Lecture Notes in Math., No. 1200,    Springer Verlag, 1986.
\bibitem[MiT]{} {\sc Milman, V.~D.} \& {\sc Tomczak-Jaegermann, N.},
  Asymptotic $l_p$ spaces and bounded distortions,
  in ``Banach Spaces'', {\sl Contemp. Math.},  144 (1993), 173--196.
\bibitem[OS.1]{}{\sc Odell, E.} \& {\sc Schlumprecht, T.},
  A Banach space block finitely universal for monotone bases, in  
  preparation.
\bibitem[OS.2]{}{\sc Odell, E.} \& {\sc Schlumprecht, T.},
  The distortion problem, to appear.
\bibitem[P.1]{}{\sc Pisier, G.},
  ``The Volume of Convex Bodies and Banach Space Geometry'',
  Cambridge Tracts in Math., 94, Cambridge Univ.~Press, 1989.
\bibitem[P.2]{}{\sc Pisier, G.},
  ``Factorization of Linear Operators and Geometry of 
  Banach Spaces'', CBMS No.~60, AMS, 1986.
\bibitem[T]{}{\sc Tomczak-Jaegermann, N.},
  ``Banach--Mazur Distances and Finite Dimensional Operator Ideals''
  Pitman Monographs, 38, Longman Scientific \& Technical, 1989.
\end{thebibliography}
\end{document}